\documentclass[12pt,a4paper,twoside]{article}

\usepackage[active]{srcltx}
\usepackage[utf8]{inputenc}
\usepackage{hyperref}
\usepackage[english]{babel}
\usepackage{authblk}

\usepackage[a4paper,hcentering,vcentering]{geometry}
\usepackage{url,epsfig,amssymb,amsthm,amsmath}

\numberwithin{equation}{section}

\newtheorem{theorems}{Theorem}
\numberwithin{theorems}{section}

\numberwithin{corollarys}{section}

\newtheorem{definition}{Definition}
\numberwithin{definition}{section}

\newtheorem{lemma}{Lemma}
\numberwithin{lemma}{section}
\newtheorem{proposition}{Proposition}
\numberwithin{proposition}{section}
\newtheorem{theoremx}{Theorem}

\newtheorem*{definitionnull}{Definition}

\newcommand{\N}{\ensuremath{\mathbb{N}}}

\newcommand{\R}{\ensuremath{\mathbb{R}}}
\newcommand{\T}{\ensuremath{\mathbb{T}}}
\newcommand{\Z}{\ensuremath{\mathbb{Z}}}

\begin{document}

\title{\LARGE{\textbf{Biasymptotically quasiperiodic solutions for time-dependent Hamiltonians}}}%
\author{Donato Scarcella}
\affil{Université Paris-Dauphine - Ceremade UMR 7534 Place du Maréchal De Lattre De Tassigny, 75016 PARIS.}
\date{}%
\maketitle


\begin{abstract}
In a previous work~\cite{Sca22a}, we consider time-dependent perturbations of a Hamiltonian vector field having an invariant torus supporting quasiperiodic solutions. Assuming the perturbation decays polynomially fast in time (when $t \to +\infty$), we prove the existence of an asymptotic KAM torus. An asymptotic KAM torus is a time-dependent family of embedded tori converging as time tends to infinity to the invariant torus associated with the unperturbed system. 
Now, it is quite natural to wonder when we have the existence of a biasymptotic KAM torus. That is, a continuous time-dependent family of embedded tori converging in the future (when $t \to +\infty$) and the past (when $t \to -\infty$) to suitable quasiperiodic invariant tori.

In this work, we go one step further. We analyze time-dependent perturbations of integrable and near-integrable Hamiltonians. Assuming the perturbation decays polynomially fast in time, we prove the existence of orbit converging to some quasiperiodic solutions in the future and the past.
\end{abstract}

\setcounter{tocdepth}{1}
\tableofcontents

\section{Introduction}

The KAM theory deals with results concerning the persistence of quasiperiodic solutions for nearly integrable Hamiltonian systems. The extraordinary work of Kolmogorov~\cite{Kol54} is followed by those of Arnold~\cite{Arn63a, Arn63b}, Moser~\cite{moser1962invariant} and others. 

The KAM theory has known many developments over the years. Concerning a recent work, the paper of Pöschel is remarkable~\cite{Posc01}. Taking into account the idea of~\cite{Mos67}, the author introduces the frequencies as independent parameters, giving a very subtle statement and an elegant proof of the theorem (in the case of real-analytic Hamiltonians). Concerning the finite differentiable case, the first proof is given by Moser~\cite{moser1962invariant}. Now, we know that it is enough to assume the system to be of class $C^k$ with $k>2(\tau+1)>2n$, for an integer $n \ge 2$ and a real number $\tau > n-1$.  We refer to~\cite{Posc82},~\cite{Sal04},~\cite{Bou20} and~\cite{Kou20}. 
On the other hand, a series of proofs are given by introducing an adapted implicit function theorem in a scale of Banach spaces, which replaces the iterative scheme inaugurated by Kolmogorov. We refer to the works of Zehnder~\cite{Zeh76, Zeh75}, Herman~\cite{bost1986tores}, Berti-Bolle~\cite{BB15} and Féjoz~\cite{Fe04, Fe16}.  One can see~\cite{M19} for the dissipative case.

In this work, we consider time-dependent perturbations of integrable Hamiltonians and Hamiltonians having a large (in the sense of measure) subset of invariant tori. Assuming the perturbation decays, on a suitable norm, polynomially fast in time, we prove the existence of orbits converging to some quasiperiodic solutions in the future ($t \to +\infty$) and the past ($t \to -\infty$). 

For the sake of clarity, we take a step back and consider a previous work~\cite{Sca22a}. In this paper, we generalized the results of Fortunati-Wiggins~\cite{FW14} and Canadell-de la Llave~\cite{CdlL15} in the particular case of time-dependent Hamiltonian systems. More specifically, we consider time-dependent perturbations of Hamiltonian vector fields (real-analytic or finitely differentiable) having an invariant torus supporting quasiperiodic solutions. Assuming the perturbation decays polynomially fast in time (when $t \to +\infty$), we prove the existence of an asymptotic KAM torus. That is a time-dependent family of embedded tori that converges as time tends to infinity to the quasiperiodic invariant torus associated with the unperturbed Hamiltonian system. We point out that Fortunati-Wiggins and Canadell-de la Llave demand exponentially fast convergence. The proof rests on the implicit function theorem. 

These kinds of perturbations appear in many physical problems. For example, in the case of a molecule disturbed by another molecule or by a laser pulse~\cite{kawai2007transition, BdlL11}. In another paper~\cite{Sca22b}, we analyse the example of the planar three-body problem perturbed by a given comet coming from and going back to infinity asymptotically along a hyperbolic Keplerian orbit (modelled as a time-dependent perturbation). In this case, the polynomial decay in time is too weak, and we can not apply the theorem proved in~\cite{Sca22a}. On the other hand, the time-dependent perturbation strongly modifies the dynamics at infinity, requiring the introduction of a generalization of the definition of asymptotic KAM torus and the formulation of another abstract theorem, whose proof is based on a Nash-Moser theorem. 

The case of time-dependent perturbations of a Hamiltonian system having an invariant torus with arbitrary dynamics is treated in~\cite{Sca22c}. Similarly to Fortunati-Wiggins and Canadell-de la Llave, we need exponential decay in time for the perturbation.

In order to state our main results, we need to introduce the following definitions.  Let $B \subset \R^n$ be a ball centred at the origin. Given $\sigma \ge 0$ and a positive integer $k \ge 0$, we consider time-dependent vector fields $X^t$ and $X^t_0$ of class $C^{\sigma + k}$ on $\T^n \times B$, for all $t \in [0, +\infty)$, and an embedding $\varphi_0$ from $\T^n$ to $\T^n \times B$ of class $C^\sigma$ such that
\begin{align}
\label{hyp1asymKAM}
& \displaystyle \lim_{t \to +\infty} |X^t - X^t_0|_{C^{\sigma + k}} = 0,\\
\label{hyp2asymKAM}
& X_0 (\varphi_0(q), t) =\partial_q \varphi_0(q)\omega \hspace{2mm} \mbox{for all $(q, t) \in \T^n \times [0, +\infty)$,}
\end{align}
where $\omega \in \R^n$ and $|\cdot|_{C^\sigma}$ is the Hölder norm. In words, the time-dependent vector field $X^t$ converges in time to the vector field $X_0^t$ that has an invariant torus supporting quasiperiodic dynamics with frequency vector $\omega$.   
\begin{definition}
\label{asymKAMtoriHol}
We assume that $(X, X_0, \varphi_0)$ satisfy~\eqref{hyp1asymKAM} and~\eqref{hyp2asymKAM}. A family of $C^\sigma$ embeddings $\varphi^t: \T^n \to \T^n \times B$ is a (positive) $C^\sigma$-asymptotic KAM torus associated to $(X, X_0, \varphi_0)$ if 
\begin{align}
\label{Hyp1}
&  \lim_{t \to +\infty}   |\varphi^t - \varphi_0|_{C^\sigma} = 0,\\
\label{Hyp2}
&  X (\varphi(q,t), t) = \partial_q \varphi(q, t) \omega + \partial_t \varphi(q, t), 
\end{align}
for all $(q, t) \in \T^n \times [0, +\infty)$. Moreover, we say that $\varphi^t$ is Lagrangian if $\varphi^t(\T^n)$ is Lagrangian for all $t$.  
\end{definition}

Obviously, if we replace $[0, +\infty)$ with $(-\infty, 0]$ and~\eqref{hyp1asymKAM},~\eqref{Hyp1} with
\begin{equation*}
 \displaystyle \lim_{t \to -\infty} |X^t - X^t_0|_{C^{\sigma + k}} = 0, \quad  \lim_{t \to -\infty}   |\varphi^t - \varphi_0|_{C^\sigma} = 0,
\end{equation*}
we obtain the definition of (negative) $C^\sigma$-asymptotic KAM torus associated to $(X, X_0, \varphi_0)$.

The above definition is due to M. Canadell and R. de la Llave~\cite{CdlL15}. In this work, they use the notion of non-autonomous KAM torus. Here, we prefer asymptotic KAM torus to point out the asymptotic properties of these objects.

As mentioned before, we proved the existence of a (positive) $C^\sigma$-asymptotic KAM torus for a time-dependent Hamiltonian vector field converging polynomially fast in time to a Hamiltonian vector field having an invariant torus supporting quasiperiodic solutions~\cite{Sca22a}. Now, it is quite natural to wonder when we have the existence of a continuous time-dependent family of embeddings $\varphi^t$ converging in the future (when $t \to +\infty$) and the past (when $t \to -\infty$) to suitable quasiperiodic invariant tori. To this end, let us introduce the definition of $C^\sigma$-biasymptotic KAM torus.

Given $\sigma \ge 0$ and a positive integer $k \ge 0$, we consider time-dependent vector fields $X^t$, $X^t_{0,+}$, $X^t_{0,-}$ of class $C^{\sigma + k}$ on $\T^n \times B$, for all $t \in \R$, and embeddings $\varphi_{0,+}$, $\varphi_{0,-}$ from $\T^n$ to $\T^n \times B$ of class $C^\sigma$ such that  
\begin{align}
\label{ToriBiasHyp1Intro}
& \displaystyle \lim_{t \to \pm\infty}  |X^t - X^t_{0,\pm}|_{C^{\sigma +k}} = 0,\\
\label{ToriBiasHyp2Intro}
& X_{0,\pm}(\varphi_{0,\pm}(q), t) = \partial_q \varphi_{0,\pm}(q)\omega_\pm \hspace{2mm} \mbox{for all $(q, t) \in \T^n \times \R$},
\end{align}
where $\omega_+$, $\omega_- \in \R^n$. 
\begin{definition}
We assume that $(X, X_{0, \pm}, \varphi_{0, \pm})$ satisfy~\eqref{ToriBiasHyp1Intro} and~\eqref{ToriBiasHyp2Intro}. For all $t \in \R$, a continuous family of $C^\sigma$ embeddings $\varphi^t : \T^n \to \T^n \times B$ is a $C^\sigma$-biasymptotic KAM torus associated to $(X, X_{0, \pm}, \varphi_{0, \pm})$ if
\begin{align*}
&  \lim_{t \to \pm\infty}   |\varphi^t - \varphi_{0, \pm}|_{C^\sigma} = 0,\\
&  X (\varphi(q,t), t) = \partial_q \varphi(q, t) \omega_+ + \partial_t \varphi(q, t), \hspace{2mm} \mbox{for all $(q,t) \in \T^n \times (0, +\infty)$}\\
&  X (\varphi(q,t), t) = \partial_q \varphi(q, t) \omega_- + \partial_t \varphi(q, t), \hspace{2mm} \mbox{for all $(q,t) \in \T^n \times (-\infty, 0)$}.
\end{align*}
Moreover, we say that $\varphi^t$ is Lagrangian if $\varphi^t(\T^n)$ is Lagrangian for all $t$.  
\end{definition}

Unfortunately, we are not able to prove the existence of $C^\sigma$-biasymptotic KAM tori without asking $\omega_+=\omega_-=\omega$ and the Hamiltonian $H:\T^n \times B \times \R \to \R$ to satisfy very strong symmetries. More specifically, we need to assume that 
\begin{equation*}
H(q + \omega t, p, t) = -H(q - \omega t, p, -t)
\end{equation*}
for all $(q, p, t)\in \T^n \times B \times \R$ and some $\omega \in \R^n$. Then, we can prove the existence of a $C^\sigma$-biasymptotic KAM torus with $\omega_+ = \omega_- = \omega$. 

In this paper, we prove weaker results. For this purpose, we introduce the following definition. 

\begin{definition}
\label{biasimsolDfifffini}
We assume that $(X, X_{0, \pm}, \varphi_{0, \pm})$ satisfy~\eqref{ToriBiasHyp1Intro} and~\eqref{ToriBiasHyp2Intro}.
An integral curve $g(t)$ of $X$ is a biasymptotically quasiperiodic solution associated to $(X, X_{0, \pm} \varphi_{0, \pm})$ if there exist $q_-$, $q_+ \in \T^n$ in such a way that 
\begin{equation}
\displaystyle \lim_{t \to \pm\infty}|g(t) - \varphi_{0, \pm}(q_\pm + \omega_\pm t)| = 0.
\end{equation} 
\end{definition}

Roughly speaking, a biasymptotically quasiperiodic solution $g$ is an orbit of the time-dependent vector field $X$, which converges to a quasiperiodic motion of frequency vector $\omega_+ \in \R^n$ in the future and to a quasiperiodic motion of frequency vector $\omega_- \in \R^n$ in the past. 

As mentioned before, this paper considers time-dependent perturbations of integrable or near-integrable Hamiltonians. Assuming the perturbation decays polynomially fast in time, we prove the existence of a large (in the sense that we will specify later) set of initial points giving rise to biasymptotically quasiperiodic solutions (see Theorem \ref{ThmInt} and Theorem \ref{ThmNearInt}). 
In both cases (integrable and near-integrable unperturbed Hamiltonian), the idea of the proof is similar.  Using different versions of the theorem in~\cite{Sca22a}, we prove the existence of a very large family of positive (resp. negative) asymptotic KAM tori. Then, we look at the intersection between these two families when $t=0$. Under suitable hypotheses on the Hamiltonian's regularity and the perturbation's smallness, it is a large set, and each point gives rise to biasymptotically quasiperiodic solutions. Contrary to our previous work~\cite{Sca22a}, here we need to assume a stronger decay in time (but always polynomial) and a smallness assumption on the perturbation.

\section{Main results}

We consider two different cases. First, we deal with time-dependent perturbations of integrable Hamiltonians, then time-dependent perturbations of autonomous Hamiltonians having a large (in the sense of measure) subset of invariant tori. 
In the integrable case, we prove the existence of biasymptotically quasiperiodic solutions for every initial condition. Concerning the other case, we show the existence of a large subset of initial conditions giving rise to biasymptotically quasiperiodic solutions.

We introduce the following notation. For every function $f$ defined on $\T^n \times B \times \R$  and, for fixed $t \in \R$, we let $f^t$ be the function defined on $\T^n \times B$ in such a way that
\begin{equation*}
f^t(q,p) = f(q,p,t).
\end{equation*}
In addition, for all fixed $p_0 \in B$, we let $f_{p_0}$ be the function defined on $\T^n \times \R$ such that
\begin{equation*}
 f_{p_0}(q,t) = f(q,p_0, t).
\end{equation*}
Obviously, for all fixed $(p,t) \in B\times \R$, we consider $f_p^t$ as the function defined on $\T^n$ such that 
\begin{equation*}
f_p^t(q) = f(q,p,t)
\end{equation*}
for all $q \in \T^n$. 

\subsection{Integrable case}

Given a positive real parameter $\sigma \ge 1$, we introduce the following space of functions.
\begin{definition}
\label{BBM}
Let $\mathcal{B}_{\sigma}$ be the space of functions $f$ defined on $\T^n \times B \times \R$ such that $f$, $\partial_p f \in C(\T^n \times B \times \R)$ and $f_p^t \in C^\sigma(\T^n)$ for all $(p,t) \in B \times \R$.
\end{definition}
We point out that $\partial_p f$ stands for partial derivatives with respect to the variables $p=(p_1,...,p_n)$ of $f$.
For all $f \in \mathcal{B}_{\sigma}$ and a real parameter $l \ge1$, we define the following special norms
\begin{align}
\label{normBMInt1}
&|f|_{\sigma , l} = \sup_{(p,t) \in B \times \R}|f^t_p|_{C^\sigma}(1 + |t|^l) + \sup_{(p,t) \in B \times \R}|\left(\partial_p f\right)^t_p|_{C^0}(1 + |t|^{l-1}),\\
 \label{normBMInt2}
&|f|_{\sigma , 0} = \sup_{(p,t) \in B \times \R}|f^t_p|_{C^\sigma}+ \sup_{(p,t) \in B \times \R}|\left(\partial_p f\right)^t_p|_{C^0}.
\end{align}
We refer to Section \ref{FS} for a series of properties of the previous norms. We need to introduce another space of functions.
\begin{definition}
\label{BBM2}
Given $\sigma \ge 1$ and an integer $k \ge 0$, we define $\mathcal{\bar B}_{\sigma, k}$ the space of functions $f$ such that 
\begin{equation*}
f \in \mathcal{B}_{\sigma +k}, \hspace{2mm} \mbox{and} \hspace{2mm} \partial^i_q f \in \mathcal{B}_{\sigma + k -i}
\end{equation*}
for all $0 \le i \le k$. 
\end{definition} 
In the latter definition, $\partial^i_q f$ stands for partial derivatives of order $i$ with respect to the variables $q=(q_1,...,q_n)$ of $f$. We use the convention $\partial^0_q f = f$. Obviously, we have $\mathcal{\bar B}_{\sigma, 0} = \mathcal{B}_{\sigma}$. Furthermore, for all $f \in \mathcal{\bar B}_{\sigma, k}$ and $l>1$, we consider the following norm
\begin{equation}
\label{normpazzerellaBM}
\left \|f \right \| _{\sigma, k, l} = \max_{0 \le i \le k}|\partial_q^i f|_{\sigma+k-i, l},
\end{equation} 
where $|\cdot|_{\sigma, l}$ is the norm defined by~\eqref{normBMInt1}

In the previous definition and Definition \ref{BBM}, $B \in \R^n$ is a ball with some unspecified radius. In what follows, we will pay attention to the radius of $B$. Let $B_r \subset \R^n$ be a ball centred at the origin with radius $r >0$. If a function defined on $\T^n \times B_r \times \R$ belongs to $\mathcal{B}_{\sigma}$, we consider that it satisfies the properties in Definition \ref{BBM} with $B$ replaced by $B_r$.

Now, we have everything we need to state the main result of this first part. 
Let $\sigma \ge 1$, $\Upsilon \ge 1$, $l>1$ and $0 < \varepsilon <1$. We consider the following Hamiltonian
\begin{equation}
\label{H1BM}
\begin{cases}
H:\T^n \times B_1 \times \R \longrightarrow \R\\
H(q,p,t) = h(p) + f(q,p,t)\\
f, \partial_p f \in \mathcal{\bar B}_{\sigma, 2},\\
|f|_{\sigma +2, 0} + \left\|\partial_q f\right\|_{\sigma,1, l+2} + \left\|\partial_p f\right\|_{\sigma,2, l+1} < \varepsilon,\\
\partial_p^2 H^t \in C^{\sigma +2}(\T^n\times B_1) \hspace{2mm} \mbox{for all fixed $t \in \R$}\\
\partial^i_{qp} \left(\partial_p^2 H\right) \in C(\T^n \times B_1 \times \R) \hspace{2mm} \mbox{for all $0 \le i \le 3$.}\\
\sup_{t \in \R}|\partial_p^2 H^t|_{C^{\sigma +2}} \le \Upsilon.
\end{cases}
\tag{$*$}
\end{equation}

For each $p\in B_1$, we consider the following trivial embedding
\begin{equation*}
\varphi_{0,p} : \T^n  \to \T^n \times B_1 , \quad \varphi_{0,p}(q) = (q,p).
\end{equation*}

\begin{theoremx}
\label{ThmInt}
Let $H$ be as in~\eqref{H1BM}. Then, there exists a time-dependent Hamiltonian $\tilde h$ such that, if $\varepsilon$ is small enough with respect to $n$, $l$, $\Upsilon$ and $|\partial_p h|_{C^1}$, for all $(q,p)\in \T^n \times B_{1 \over 2}$ there exist $p_-$, $p_+ \in B_1$ and a biasymptotically quasiperiodic solution $g(t)$ associated to $(X_H, X_{\tilde h}, \varphi_{0,p_\pm})$ such that $g(0) = (q,p)$.
\end{theoremx}

Concerning the regularity of the time-dependent perturbation $f$, if we assume that for all fixed $(p,t) \in B_1 \times \R$ 
\begin{equation*}
f^t_p, \partial_p f^t_p \in C^{\sigma+2}(\T^n) \quad \mbox{and} \quad \partial^d_q\partial^m_p f \in C(\T^n \times B \times \R)
\end{equation*}
for all $0 \le d +m \le 4$ with $0 \le d \le 2$ and $0 \le m \le 2$, then one can prove that $f, \partial_p f \in \mathcal{\bar B}_{\sigma, 2}$. Obviously, we point out that if we assume the stronger hypothesis
\begin{equation*}
f \in C^{\sigma  + 3} (\T^n \times B_1 \times \R) \hspace{2mm} \mbox{and} \hspace{2mm} \partial_p^2 H \in C^{\sigma  + 2} (\T^n \times B_1 \times \R),
\end{equation*}
then the regularity assumptions of the previous theorem are satisfied. 

Instead of proving this theorem directly, we are going to deduce it from another theorem. 
Let $\sigma \ge 1$, $\Upsilon \ge 1$, $l>1$ and $0 < \varepsilon <1$. We consider the following family of Hamiltonians
\begin{equation}
\label{H3BM}
\begin{cases}
H:\T^n \times B_{1 \over 4} \times \R \times B_{3 \over 4} \longrightarrow \R\\
H(\theta, I, t;p_0) = e(p_0) + \omega(p_0) \cdot I \\
\hspace{23mm} + a(\theta, t;p_0) + b(\theta, t;p_0) \cdot I + m(\theta, I, t;p_0) \cdot I^2\\
 \omega \in C^1(B_{3 \over 4}), \hspace{5mm}a, b \in \mathcal{\bar B}_{\sigma, 2},\\
|a|_{\sigma +2, 0} + \left\|\partial_\theta a\right\|_{\sigma,1, l+2} < \varepsilon, \hspace{2mm} \left\|b\right\|_{\sigma,2, l+1} < \varepsilon,\\
\partial_I^2 H^t \in C^{\sigma +2}(\T^n\times B_{1 \over 4} \times B_{3 \over 4}) \hspace{2mm} \mbox{for all fixed $t \in \R$}\\
\partial^i_{\theta I p_0} \left(\partial_I^2 H\right) \in C(\T^n \times B_{1 \over 4} \times \R \times B_{3 \over 4})  \hspace{2mm} \mbox{for all $0 \le i \le 3$.}\\
\sup_{t \in \R}|\partial_I^2 H^t|_{C^{\sigma +2}} \le \Upsilon.
\end{cases}
\tag{$\star$}
\end{equation} 
We define the following family of trivial embeddings
\begin{equation}
\label{trivialembeddingBM}
\psi_0 : \T^n \times B_{3 \over 4} \longrightarrow \T^n \times B_{1 \over 4}, \quad \psi_0(\theta, p_0) = (\theta,0)
\end{equation}
and we consider the following family of Hamiltonians $\tilde h:\T^n \times B_{1 \over 4} \times \R \times B_{3 \over 4} \to \R$ such that
\begin{equation*}
\tilde h(\theta, I, t;p_0) = e(p_0) + \omega(p_0) \cdot I + m(\theta, I, t;p_0) \cdot I^2.
\end{equation*}

\begin{theorems}
\label{Thm3BM}
Let $H$ be as in~\eqref{H3BM}. Then, if $\varepsilon$ is sufficiently small with respect to $n$, $l$, $\Upsilon$ and $|\omega|_{C^1}$, for all fixed $p_0 \in B_{3 \over 4}$, there exists a positive $C^\sigma$-asymptotic KAM torus $\psi^t_{+p_0}$ and a negative $C^\sigma$-asymptotic KAM torus $\psi^t_{-p_0}$ associated to $(X_{H_{p_0}}, X_{\tilde h_{p_0}}, \psi_{0,p_0})$. Moreover, letting
\begin{equation*}
\psi^t_\pm : \T^n \times B_{3 \over 4} \longrightarrow \T^n \times B_{1 \over 4},\quad \psi^t_\pm (q,p_0) = \psi^t_{\pm,p_0}(q),
\end{equation*}
there exists a constant $C_0$ depending on $n$, $l$, $\Upsilon$ and $|\omega|_{C^1}$ such that 
\begin{equation}
\label{stimethm3BM}
\sup_{t \ge 0} |\psi^t_+ - \psi_0|_{C^1} < C_0 \varepsilon, \quad \sup_{t \le 0} |\psi^t_- - \psi_0|_{C^1} < C_0 \varepsilon.
\end{equation}
\end{theorems}

For clarity, we point out that, for all fixed $p_0 \in B_{3 \over 4}$,  $\psi_{0,p_0}(q) = \psi_0(q,p_0)$ for all $q \in \T^n$.

\subsection{Near integrable case}\label{ResultsBMnotIntCas}

Let $A \subset \R^n$ be a closed set and $E$ equal to $\T^n$ or $\T^n \times B$. For every function $f : E \times A \subset \R^n \times \R^n \to \R$ we introduce the following Lipschitz norm
\begin{equation*}
|f|_{L(A)} = \sup_{z \in B}\left(\sup_{x,y \in A,  x\ne y} {\left|f(z, x) - f(z, y)\right| \over \left|x-y\right|}\right) + |f|_{C^0}.
\end{equation*}
Let a real parameter $\sigma \ge 1$ and $A \subset \R^n$. We consider a suitable space of functions.
\begin{definition}
\label{DBM}
Let $\mathcal{D}_{\sigma}$ be the space of functions $f$ defined on $\T^n \times A \times \R$ such that $f \in C(\T^n \times A \times \R)$ and $f_p^t \in C^\sigma(\T^n)$ for all $(p,t) \in A \times \R$.
\end{definition}

For all $f \in \mathcal{D}_{\sigma}$ and $l \ge 1$, we define 
\begin{align}
\label{normBMNotInt}
&|f|_{\sigma , l, L(A)} = \sup_{(p,t) \in A \times \R}|f^t_p|_{C^\sigma}(1 + |t|^l) + \sup_{t \in  \R}|f^t|_{L(A)}(1 + |t|^{l-1}),\\
&|f|_{\sigma , 0, L(A)} = \sup_{(p,t) \in A \times \R}|f^t_p|_{C^\sigma}+ \sup_{t \in  \R}|f^t|_{L(A)}.
\end{align}
Also in this case, we refer to Section \ref{FS} for a series of properties of the previous norms. As in the previous theorem, we define the following space of functions.
\begin{definition}
\label{DBM2}
Given $\sigma \ge 1$ and an integer $k \ge 0$, we define $\mathcal{\bar D}_{\sigma, k}$ the space of functions $f$ such that 
\begin{equation*}
f \in \mathcal{D}_{\sigma +k}, \hspace{2mm} \mbox{and} \hspace{2mm} \partial^i_q f \in \mathcal{D}_{\sigma + k -i}
\end{equation*}
for all $0 \le i \le k$. 
\end{definition} 
 Furthermore, for all $f \in \mathcal{\bar D}_{\sigma, k}$ and $l>1$, we consider the following norm
\begin{equation}
\label{normpazzerellaBMNotInt}
\left \|f \right \| _{\sigma, k, l, L(A)} = \max_{0 \le i \le k}|\partial_q^i f|_{\sigma+k-i, l, L(A)},
\end{equation} 
where $|\cdot|_{\sigma, l, L(A)}$ is the norm defined by~\eqref{normBMNotInt}

Here, the following Hamiltonian is defined on $\T^n \times B_r \times \R$, for some $r>0$. Then, similarly to the previous case, if a function defined on $\T^n \times B_r \times \R$ belongs to the space $\mathcal{D}_{\sigma}$, we consider that it satisfies the properties in Definition \ref{BBM} with $A$ replaced by $B_r$.

 Now, we can state the main result of this second part. Let $\sigma \ge 1$, $\Upsilon \ge 1$, $l>1$, $0 < \varepsilon <1$ and $\mu>0$. We consider the following Hamiltonian
\begin{equation}
\label{H2BM}
\begin{cases}
H:\T^n \times B_1\times \R \longrightarrow \R\\
H(q,p,t) = h(p) + R(q,p) + f(q,p,t)\\
f, \partial_p f \in \mathcal{\bar D}_{\sigma, 2},\\
|f|_{\sigma +2, 0, L(B_1)} + \left\|\partial_q f\right\|_{\sigma,1, l+2, L(B_1)} + \left\|\partial_p f\right\|_{\sigma,2, l+1, L(B_1)} < \varepsilon,\\
D \subset B_1, \quad \mathrm{Leb}(B_1\backslash D) < \mu,\\
R \in C^2(\T^n \times B_1)\\
R(q,p) = \partial_p R(q,p) =0 \hspace{2mm} \mbox{for all $(q,p) \in \T^n \times D$},\\
\partial_p^2 H^t \in C^{\sigma +2}(\T^n\times B_1) \hspace{2mm} \mbox{for all fixed $t \in \R$}\\
\partial^i_{qp} \left(\partial_p^2 H\right) \in C(\T^n \times B_1 \times \R) \hspace{2mm} \mbox{for all $0 \le i \le 2$.}\\
\sup_{t \in \R}|\partial_p^2 H^t|_{C^{\sigma +2}} \le \Upsilon,
\end{cases}
\tag{$**$}
\end{equation}

For each $p\in B_1$, we recall that $\varphi_{0,p}$ is the following trivial embedding
\begin{equation*}
\varphi_{0,p} : \T^n  \to \T^n \times B_1 , \quad \varphi_{0,p}(q) = (q,p).
\end{equation*}

\begin{theoremx}
\label{ThmNearInt}
Let $H$ be as in~\eqref{H2BM}. Then, there exists a time-dependent Hamiltonian $\tilde h$ such that, for $\varepsilon$ small enough with respect to $n$, $l$, $\Upsilon$, $|\partial_p h|_{L(D)}$ and $\mu$, we have the existence of a set $\mathcal{W} \subset \T^n \times B_1$ in such a way that, for all $(q,p) \in \mathcal{W}$, there exist $p_-$, $p_+ \in D$ and a biasymptotically quasiperiodic solution $g$ associated to $(X_H, X_{\tilde h}, \varphi_{0,p_\pm})$ such that $g(0) =(q,p)$. Moreover, 
\begin{equation*}
\mathrm{Leb} \left(\left(\T^n \times B_1\right)\backslash \mathcal{W}\right) \le 4 \mu.
\end{equation*}
\end{theoremx}

First of all, if we assume that for all fixed $(p,t) \in B_1 \times \R$ 
\begin{equation*}
f^t_p, \partial_p f^t_p \in C^{\sigma+2}(\T^n) \quad \mbox{and} \quad \partial^d_q\partial^m_p f \in C(\T^n \times B \times \R)
\end{equation*}
for all $0 \le d +m \le 3$ with $0 \le d \le 2$ and $0 \le m \le 1$, then one can prove that $f, \partial_p f \in \mathcal{\bar D}_{\sigma, 2}$.

Concerning the hypothesis on the autonomous part $h+R$ of the previous Hamiltonian $H$, we point out that it is not artificial. Pöschel, in his work~\cite{Posc82}, considers a $C^\infty$ small perturbation $H_1$ of a real-analytic integrable non-degenerate Hamiltonian $H_0$ of the form
\begin{equation*}
H:\T^n \times B \longrightarrow \R, \quad H(q,p) = H_0(p) + H_1(q,p).
\end{equation*}
The author proves the existence of a $C^\infty$-symplectomorphism $\phi$ such that 
\begin{equation*}
H \circ \phi:\T^n \times B \longrightarrow \R, \quad H(q,p) = h(p) + R(q,p),
\end{equation*}
where $h$ is close to $H_0$ and the infinite jet of $R$ vanishes along $\T^n \times D_\gamma$. Moreover, $D_ \gamma$ is a suitable subset of $B$ such that 
\begin{equation*}
\mathrm{Leb}(B\backslash D_\gamma) \le C\gamma^2,
\end{equation*}
for an appropriate constant $C$ and a positive parameter $\gamma >0$. By referring to this work of Pöschel, in Theorem \ref{ThmNearInt}, we take $D=D_\gamma$ and $\mu = C\gamma^2$.

Also in this case,  we are going to deduce the latter from another theorem. Let $\sigma \ge 1$, $\Upsilon \ge 1$, $l>1$, $0 < \varepsilon <1$, $\mu >0$ and $0<\delta<1$ with $\delta \le \mu$. We consider the following family of Hamiltonians
\begin{equation}
\label{H4BM}
\begin{cases}
H:\T^n \times B_{\delta} \times \R \times D' \longrightarrow \R\\
H(\theta, I, t;p_0) = e(p_0) + \omega(p_0) \cdot I \\
\hspace{23mm} + a(\theta, t;p_0) + b(\theta, t;p_0) \cdot I + m(\theta, I, t;p_0) \cdot I^2\\
 \omega \in C(D'), \quad |\omega|_{L(D')}<\infty\\
a, b \in \mathcal{\bar D}_{\sigma, 2},\\
|a|_{\sigma +2, 0, L(D')} + \left\|\partial_\theta a\right\|_{\sigma,1, l+2, L(D')} < \varepsilon, \hspace{2mm} \left\|b\right\|_{\sigma,2, l+1, L(D')} < \varepsilon,\\
\partial_I^2 H^t_{p_0} \in C^{\sigma +2}(\T^n\times B_\delta) \hspace{2mm} \mbox{for all fixed $(t, p_0) \in \R \times D'$}\\
\partial^i_{\theta I } \left(\partial_I^2 H\right) \in C(\T^n \times B_\delta \times \R \times D')  \hspace{2mm} \mbox{for all $0 \le i \le 2$.}\\
\sup_{(t, p_0) \in \R \times D'}|\partial_I^2 H_{p_0}^t|_{C^{\sigma +2}} \le \Upsilon,\\
\sup_{0 \le i \le 2}\left(\sup_{t \in \R} \left|\partial^i_{\theta I}\left(\partial_I^2 H^t\right)\right|_{L(D')}\right) \le \Upsilon.
\end{cases}
\tag{$\star\star$}
\end{equation} 
We define the following family of trivial embeddings
\begin{equation*}
\psi_0 : \T^n \times D' \longrightarrow \T^n \times B_\delta, \quad \psi_0(\theta, p_0) = (\theta,0)
\end{equation*}
and we consider the following family Hamiltonians $\tilde h:\T^n \times B_{\delta} \times \R \times D' \to \R$ such that
\begin{equation*}
\tilde h(\theta, I, t;p_0) = e(p_0) + \omega(p_0) \cdot I + m(\theta, I, t;p_0) \cdot I^2.
\end{equation*}

\begin{theorems}
\label{Thm4BM}
Let $H$ be as in~\eqref{H4BM}. Then, if $\varepsilon$ is sufficiently small with respect to $n$, $l$, $\Upsilon$ and $|\omega|_{L(D')}$, for all fixed $p_0 \in D'$, there exists a positive $C^\sigma$-asymptotic KAM torus $\psi^t_{+p_0}$ and a negative $C^\sigma$-asymptotic KAM torus $\psi^t_{-p_0}$ associated to $(X_{H_{p_0}}, X_{\tilde h_{p_0}}, \psi_{0,p_0})$. Moreover, letting
\begin{equation*}
\psi^t_\pm : \T^n \times D' \longrightarrow \T^n \times B_\delta,\quad \psi^t_\pm (q,p_0) = \psi^t_{\pm,p_0}(q),
\end{equation*}
there exists a constant $C_0$ depending on $n$, $l$, $\Upsilon$ and $|\omega|_{L(D')}$ such that 
\begin{equation}
\label{stimethm4BM}
\sup_{t \ge 0} |\psi^t_+ - \psi_0|_{L(\T^n \times D')} < C_0 \varepsilon, \quad \sup_{t \le 0} |\psi^t_- - \psi_0|_{L(\T^n \times D')} < C_0 \varepsilon.
\end{equation}
\end{theorems}

We point out that Theorem \ref{Thm4BM} and Theorem \ref{Thm3BM} are different versions of the result proved in~\cite{Sca22a}. Here, we need to require more conditions to obtain more information. More specifically, we assume a stronger decay in time and a smallness assumption on the perturbative terms compared to our previous work. The reason is the control given by~\eqref{stimethm3BM} and~\eqref{stimethm4BM}, which is the key to proving the main theorems of this work (Theorem \ref{ThmInt} and Theorem \ref{ThmNearInt}). 

\section{Functional Setting}\label{FS}

In what follows, we provide a series of properties of the norms introduced in the previous section (see~\eqref{normBMInt1} and~\eqref{normBMNotInt}). But, first, we need to recall the notion of Hölder classes of functions $C^\sigma$ and some properties. 

We consider $G$ as an open subset of $\R^n$. Let $k \ge 0$ be a positive integer, we define $C^k(G)$ as the spaces of functions $f: G \to \R$ with continuous partial derivatives $\partial^\alpha f \in C^0(G)$ for all $\alpha \in \N^n$ with $|\alpha|=\alpha_1+...+\alpha_n \le k$. Then, for all $f \in C^k(G)$, we have the following norm
\begin{equation*}
|f|_{C^k} = \sup_{|\alpha|\le k}|\partial^\alpha f|_{C^0},
\end{equation*}
where $|\partial^\alpha f|_{C^0} = \sup_{x \in G}|\partial^\alpha f(x)|$ denotes the sup norm. Given $\sigma=k+\mu$, with $k \in \Z$, $k \ge 0$ and $0 < \mu <1$, we define the Hölder spaces $C^\sigma(G)$ as the spaces of functions $f\in C^k(G)$ verifying
\begin{equation}
\label{Holdernorm}
|f|_{C^\sigma} = \sup_{|\alpha|\le k}|\partial^\alpha f|_{C^0} + \sup_{|\alpha| = k}{|\partial^\alpha f(x) - \partial^\alpha f(y)| \over |x-y|^\mu}<\infty.
\end{equation}
We have the following proposition.  In this work, $C(\cdot)$ stands for constants depending on $n$ and the other parameters in brackets.  On the other hand, $C$ denotes constants depending only on $n$.

\begin{proposition}
\label{Holder}
We consider $f$, $g \in C^\sigma(G)$ and $\sigma \ge 0$.
\begin{enumerate}
\item For all $\beta \in \N^{n}$, if $|\beta| + s = \sigma$ then  $\left|{\partial^{|\beta|} \over \partial{x_1}^{\beta_1}... \partial{x_n}^{\beta_n}} f \right|_{C^s} \le |f|_{C^\sigma}$.\\
\item  $|fg|_{C^\sigma} \le C(\sigma)\left(|f|_{C^0}|g|_{C^\sigma} + |f|_{C^\sigma}|g|_{C^0}\right)$. 
\end{enumerate}
Concerning composite functions. Let $z$ be defined on $G_1 \subset \R^n$ and takes its values on $G_2 \subset \R^{n}$ where $f$ is defined. 
If $\sigma \ge 1$ and $f \in C^\sigma (G_2)$, $z \in C^\sigma (G_1)$ then $f\circ z \in C^\sigma(G_1)$ 
\begin{enumerate}
\item[3.] $|f \circ z|_{C^\sigma} \le C(\sigma) \left(|f|_{C^\sigma}|z|^\sigma_{C^1} + |f|_{C^1}| z|_{C^\sigma}+ |f|_{C^0}\right)$.
\end{enumerate}
\end{proposition}
\begin{proof}
We refer to~\cite{Hor76} for the proof.
\end{proof}

We recall that $B \subset \R^n$ is a ball centred at the origin, and the space of functions $\mathcal{B}_{\sigma}$ is introduced in Definition \ref{BBM}.  The norm defined by~\eqref{normBMInt1} satisfies the following properties.

\begin{proposition}
\label{normpropertiesBM}
Given $\sigma \ge 1$, for all $f$, $g \in \mathcal{B}_{\sigma}$ and positive $l$, $m \ge 1$
\begin{enumerate}
\item[a.] $|f|_{\sigma, l} \le |f|_{s, l}$ for all $1 \le \sigma \le s$, 
\item[b.] $|f|_{\sigma, l} \le C(l,m) |f|_{\sigma, l+m}$ \\
\item[c.]  $|fg|_{\sigma, l+m} \le C(\sigma)\left(|f|_{0,l}|g|_{\sigma,m} + |f|_{\sigma,l}|g|_{0,m}\right)$. 
\end{enumerate}
Moreover, we consider $\tilde g\in \mathcal{B}_{\sigma}$ such that, for all $(q,p,t) \in \T^n \times B \times \R$, $\tilde g(q,p,t) = (g(q,p,t), p,t)$.  Then $f \circ \tilde g \in \mathcal{B}_{\sigma}$ and  
\begin{enumerate}
\item[d.] $|f \circ \tilde g|_{\sigma, l+m} \le C(\sigma) \left(|f|_{\sigma,l} |g|^\sigma_{1,m} + |f|_{1,l}|g|_{\sigma,m} +  |f|_{0, l+m}  \right)$.
\end{enumerate}
\end{proposition}
Before the proof, we observe that the previous properties are still verified when $l =m =0$ or only one of the two parameters, $l$ and $m$, is zero.
\begin{proof}
The proof rests on Proposition \ref{Holder}. Properties \textit{a.} and \textit{b.} are obvious. Then, we verify the others.

\textit{c}. For all fixed $(p,t) \in B \times \R$, by property \textit{2.} of Proposition \ref{Holder}
\begin{eqnarray*}
 \left|f^t_p g^t_p \right|_{C^\sigma} \left(1 + |t|^{l+m}\right) &\le&  C(\sigma)\left(|f^t_p|_{C^0}|g^t_p|_{C^\sigma} + |f^t_p|_{C^\sigma}|g^t_p|_{C^0}\right)\left(1 + |t|^l \right)\left(1 + |t|^m \right) \\
&\le& C(\sigma)\Big(|f_p^t|_{C^0}\left(1 + |t|^l \right) |g_p^t|_{C^\sigma}\left(1 + |t|^m \right)\\
&+& |f_p^t|_{C^\sigma}\left(1 + |t|^l \right)|g_p^t|_{C^0}\left(1 + |t|^m \right) \Big)\\
&\le& C(\sigma)\left(|f|_{0,l}|g|_{\sigma,m} + |f|_{\sigma,l}|g|_{0,m}\right)
\end{eqnarray*}
where in the second line we use $ \left(1 + |t|^{l+m} \right) \le  \left(1 + |t|^l \right) \left(1 + |t|^m \right)$.
Taking the sup for all $(p,t) \in B \times \R$ on the left-hand side of the latter, we obtain 
\begin{equation*}
 \sup_{(p,t) \in B \times \R} \left|f^t_p g^t_p \right|_{C^\sigma} \left(1 + |t|^{l+m}\right) \le  C(\sigma)\left(|f|_{0,l}|g|_{\sigma,m} + |f|_{\sigma,l}|g|_{0,m}\right).
\end{equation*}
It remains to prove that the second term of the norm (see the right-hand side of~\eqref{normBMInt1}) also satisfies the same estimate. For all fixed $(p,t) \in B \times \R$
\begin{eqnarray*}
 \left|\left(\partial_p \left(f g\right)\right)_p^t \right|_{C^0} \left(1 + |t|^{l+m-1}\right) &=&   \left|\left(\partial_pf\right)^t_p g^t_p + f^t_p \left(\partial_p g\right)_p^t \right|_{C^0} \left(1 + |t|^{l+m-1}\right) \\
&\le& \left(\left|\left(\partial_pf\right)^t_p g^t_p\right|_{C^0}  + \left|f^t_p \left(\partial_p g\right)_p^t \right|_{C^0} \right)\left(1 + |t|^{l+m-1}\right)\\
&\le& C\left|\left(\partial_pf\right)^t_p \right|_{C^0}\left(1 + |t|^{l-1}\right) \left|g^t_p\right|_{C^0}\left(1 + |t|^{m}\right)\\
&+& C\left|\left(\partial_p g\right)_p^t \right|_{C^0}\left(1 + |t|^{m-1}\right) \left|f^t_p\right|_{C^0} \left(1 + |t|^l\right) \\
&\le& C\left(|f|_{0,l}|g|_{\sigma,m} + |f|_{\sigma,l}|g|_{0,m}\right)
\end{eqnarray*}
and taking the sup for all $(p,t) \in B \times \R$ on the left-hand side of the latter, we verify the estimate. 

\textit{d}. For all fixed $(p,t) \in B \times \R$ and thanks to property \textit{3.} of Proposition \ref{Holder}
\begin{eqnarray*}
 |f^t_p\circ g_p^t|_{C^\sigma}\left(1 + |t|^{l+m} \right)  &\le & C(\sigma) \left (\left|f_p^t\right|_{C^\sigma} \left|g_p^t\right|^\sigma_{C^1} + \left|f_p^t\right|_{C^1} \left|g_p^t\right|_{C^\sigma} + \left|f_p^t\right|_{C^0}\right) \left(1 + |t|^{l+m} \right)  \\
&\le& C(\sigma) \Big(\left|f_p^t\right|_{C^\sigma}\left(1 + |t|^l \right)  \left|g_p^t\right|^\sigma_{C^1}\left(1 + |t|^m \right)^\sigma \\
 &+& \left|f_p^t\right|_{C^1}\left(1 + |t|^l \right) \left|g_p^t\right|_{C^\sigma}\left(1 + |t|^m \right)\\
 &+& \left|f_p^t\right|_{C^0}\left(1 + |t|^{l+m} \right) \Big) \\
&\le& C(\sigma) \left(|f|_{\sigma,l} |g|^\sigma_{1,m} + |f|_{1,l}|g|_{\sigma,m} +  |f|_{0, l+m}  \right)
\end{eqnarray*}
where in the second line we use $\left(1 + |t|^m \right) \le \left(1 + |t|^m \right)^\sigma$. Taking the sup for all $(p,t) \in B_1 \times \R$ on the left-hand side of the latter,
\begin{equation*}
\sup_{(p,t) \in B \times \R}|f^t_p\circ g_p^t|_{C^\sigma}\left(1 + |t|^{l+m} \right) \le C(\sigma) \left(|f|_{\sigma,l} |g|^\sigma_{1,m} + |f|_{1,l}|g|_{\sigma,m} +  |f|_{0, l+m}  \right). 
\end{equation*}
Concerning the second term of the norm (see~\eqref{normBMInt1}), for all fixed $(p,t) \in B \times \R$

\begin{eqnarray*}
 \left|\left(\partial_p \left(f\circ \tilde  g\right)\right)_p^t\right|_{C^0}\left(1 + |t|^{l+m-1} \right)  &=& \left |\left(\partial_q f\right)^t_p\circ  g^t_p \left(\partial_p g\right)^t_p\right|_{C^0}\left(1 + |t|^{l+m-1} \right)\\
&+& \left |\left(\partial_p f\right)^t_p\circ  g^t_p \right|_{C^0}\left(1 + |t|^{l+m-1} \right) \\
&\le& C \left |\left(\partial_q f\right)^t_p\right|_{C^0}\left(1 + |t|^l \right)\left| \left(\partial_p g\right)^t_p\right|_{C^0}\left(1 + |t|^{m-1} \right)\\
&+& \left |\left(\partial_p f\right)^t_p \right|_{C^0}\left(1 + |t|^{l+m-1} \right)\\
&\le&C\left( |f|_{1,l} |g|_{\sigma, m} + |f|_{0, l+m}\right).
\end{eqnarray*}
Taking the sup for all $(p,t) \in B \times \R$ on the left-hand side of the above inequality, 
\begin{eqnarray*}
\sup_{(p,t) \in B \times \R}\left|\left(\partial_p \left(f\circ \tilde  g\right)\right)_p^t\right|_{C^0}\left(1 + |t|^{l+m-1} \right)  &\le& C\left( |f|_{1,l} |g|_{\sigma, m} + |f|_{0, l+m}\right).
\end{eqnarray*}
\end{proof}

As for the norm~\eqref{normBMInt1}, the following proposition contains several properties of the norm defined by~\eqref{normBMNotInt}. First, we recall that $A \subset \R^n$, and $\mathcal{D}_{\sigma}$ is the space of function of Definition \ref{DBM}.

\begin{proposition}
\label{normpropertiesBMNotInt}
Given $\sigma \ge 1$, for all $f$, $g \in \mathcal{D}_{\sigma}$ and positive $l$, $m \ge 1$
\begin{enumerate}
\item[a.] $|f|_{\sigma, l, L(A)} \le |f|_{s, l, L(A)}$ for all $1 \le \sigma \le s$, 
\item[b.] $|f|_{\sigma, l, L(A)} \le C(l,m) |f|_{\sigma, l+m, L(A)}$ \\
\item[c.]  $|fg|_{\sigma, l+m, L(A)} \le C(\sigma)\left(|f|_{0,l, L(A)}|g|_{\sigma,m, L(A)} + |f|_{\sigma,l, L(A)}|g|_{0,m, L(A)}\right)$. 
\end{enumerate}
Moreover, we consider $\tilde g\in \mathcal{D}_{\sigma}$ such that, for all $(q,p,t) \in \T^n \times A \times \R$, $\tilde g(q,p,t) = (g(q,p,t), p,t)$.  Then $f \circ \tilde g \in \mathcal{D}_{\sigma}$ and  
\begin{enumerate}
\item[d.] $|f \circ \tilde g|_{\sigma, l+m, L(A)} \le C(\sigma) \left(|f|_{\sigma,l, L(A)} |g|^\sigma_{1,m, L(A)} + |f|_{1,l, L(A)}|g|_{\sigma,m, L(A)} +  |f|_{0, l+m, L(A)}  \right)$.
\end{enumerate}
\end{proposition}
The previous properties are still verified when $l =m =0$ or only one of the two parameters, $l$ or $m$, is zero.
\begin{proof}
The proof is quite similar to that of Proposition \ref{normpropertiesBM}. The first two properties \textit{a.} and \textit{b.} are obvious. Hence, we prove \textit{c.} and \textit{d.}

\textit{c}. Similarly to Proposition \ref{normpropertiesBM}, thanks to property \textit{2.} of Proposition \ref{Holder}, one has
\begin{equation*}
 \sup_{(p,t) \in A \times \R} \left|f^t_p g^t_p \right|_{C^\sigma} \left(1 + |t|^{l+m}\right) \le  C(\sigma)\left(|f|_{0,l, L(A)}|g|_{\sigma,m, L(A)} + |f|_{\sigma,l, L(A)}|g|_{0,m, L(A)}\right).
\end{equation*}
We want to prove that the same estimate is also verified for the second term of the norm (see the right-hand side of~\eqref{normBMNotInt}). For all fixed $t \in \R$ and $x$, $y \in A$ such that $x \ne y$
\begin{align*}
&{|f^t(q, x) g^t(q, x) - f^t(q,y)g^t(q, y)| \over |x-y|}(1 + |t|^{l +m -1})\\
&={|f^t(q, x) g^t(q, x)-f^t(q, y) g^t(q, x)+ f^t(q, y) g^t(q, x) - f^t(q,y)g^t(q, y)| \over |x-y|}(1 + |t|^{l +m -1})\\
&\le C\left(|f^t|_{L(A)}\left(1 + |t|^{l-1}\right) |g^t_x|_{C^0}\left(1 + |t|^m\right) +  |g^t|_{L(A)}\left(1 + |t|^{m-1}\right) |f^t_y|_{C^0}\left(1 + |t|^l\right) \right)\\
&\le  C\left(|f|_{0,l, L(A)}|g|_{\sigma,m, L(A)} + |f|_{\sigma,l, L(A)}|g|_{0,m, L(A)}\right).
\end{align*}
Taking the sup for all $q \in \T^n$ and $x$, $y \in A$ with $x \ne y$ on the left-hand side of the latter and then for all $t \in \R$, we prove \textit{c.}

\textit{d}. Also in this case, similarly to Proposition \ref{normpropertiesBM} and by property \textit{3.} of Proposition \ref{Holder}, one has
\begin{eqnarray*}
\sup_{(p,t) \in A \times \R}|f^t_p\circ g_p^t|_{C^\sigma}\left(1 + |t|^{l+m} \right) &\le& C(\sigma)|f|_{\sigma,l, L(A)} |g|^\sigma_{1,m, L(A)} \\
&+& C(\sigma)\left(|f|_{1,l, L(A)}|g|_{\sigma,m, L(A)} +  |f|_{0, l+m, L(A)}  \right). 
\end{eqnarray*}
Now, we estimate the second term of the norm (see~\eqref{normBMNotInt}). For all fixed $t \in \R$ and $x$, $y \in A$ such that $x \ne y$
\begin{align*}
&{|f^t(g^t(q, x), x) - f^t(g^t(q, y), y)| \over |x-y|}(1 + |t|^{l +m -1})\\
&={|f^t(g^t(q, x), x) - f^t(g^t(q, x), y) + f^t(g^t(q, x), y) - f^t(g^t(q, y), y)| \over |x-y|}(1 + |t|^{l +m -1})\\
&\le{|f^t(g^t(q, x), x) - f^t(g^t(q, x), y)| \over |x-y|}(1 + |t|^{l +m -1})\\
&+{|f^t(g^t(q, x), y) - f^t(g^t(q, y), y)| \over |x-y|}(1 + |t|^{l +m -1})\\
&\le C\left(|f^t|_{L(A)}\left(1 + |t|^{l+m-1}\right) + |f^t_y|_{C^1}\left(1 + |t|^l\right) |g^t|_{L(A)}\left(1 + |t|^{m-1}\right)  \right)\\
&\le  C\left(|f|_{0,l+m, L(A)} + |f|_{1,l, L(A)}|g|_{\sigma,m, L(A)}\right),
\end{align*}
where we used $|f^t(g^t(q, x), y) - f^t(g^t(q, y), y)| \le |\left(\partial_q f\right)^t_y|_{C^0}|g^t|_{L(A)}|x-y|$.
Taking the sup for all $q \in \T^n$ and $x$, $y \in A$ with $x \ne y$ on the left-hand side of the latter and then for all $t \in \R$, we conclude the proof of this lemma. 
\end{proof}

\section{Biasymptotically quasiperiodic solutions}

We recall the definition of biasymptotically quasiperiodic solutions. 
Given $\sigma \ge 0$ and a positive integer $k \ge 0$, we consider time-dependent vector fields $X^t$, $X^t_{0,+}$, $X^t_{0,-}$ of class $C^{\sigma + k}$ on $\T^n \times B$, for all $t \in \R$, and embeddings $\varphi_{0,+}$, $\varphi_{0,-}$ from $\T^n$ to $\T^n \times B$ of class $C^\sigma$ such that  
\begin{align}
\label{ToriBiasHyp1}
& \displaystyle \lim_{t \to \pm\infty}  |X^t - X^t_{0,\pm}|_{C^{\sigma +k}} = 0,\\
\label{ToriBiasHyp2}
& X_{0,\pm}(\varphi_{0,\pm}(q), t) = \partial_q \varphi_{0,\pm}(q)\omega_\pm \hspace{2mm} \mbox{for all $(q, t) \in \T^n \times \R$},
\end{align}
where $\omega_+$, $\omega_- \in \R^n$. 

\begin{definitionnull}[Definition \ref{biasimsolDfifffini}]
We assume that $(X, X_{0, \pm}, \varphi_{0, \pm})$ satisfy~\eqref{ToriBiasHyp1} and~\eqref{ToriBiasHyp2}.
An integral curve $g(t)$ of $X$ is a biasymptotically quasiperiodic solution associated to $(X, X_{0, \pm} \varphi_{0, \pm})$ if there exist $q_-$, $q_+ \in \T^n$ in such a way that 
\begin{equation}
\displaystyle \lim_{t \to \pm\infty}|g(t) - \varphi_{0, \pm}(q_\pm + \omega_\pm t)| = 0.
\end{equation} 
\end{definitionnull}

This section is devoted to the proof of the following obvious proposition. Let $X$, $X_{0,+}$, $X_{0,-}$, $\varphi_{0,+}$ and $\varphi_{0,-}$ be as in the previous definition
\begin{proposition}
\label{proplink2BM}
We assume the existence of a positive $C^\sigma$-asymptotic KAM torus $\varphi^t_+$ associated to $(X, X_{0,+} \varphi_{0,+})$, a negative $C^\sigma$-asymptotic KAM torus $\varphi^t_-$ associated to $(X, X_{0,-} \varphi_{0,-})$ and $q_+$, $q_- \in \T^n$ in such a way that
\begin{equation}
\label{varphi0+-}
\varphi^0_+(q_+) = \varphi^0_-(q_-).
\end{equation}
Then, letting
\begin{equation}
\label{g}
g(t)=\begin{cases}   \varphi_+^t (q_+ + \omega_+t) \quad \mbox{for all $t \ge 0$}\\
\varphi_-^t (q_- + \omega_-t) \quad \mbox{for all $t \le 0$},
\end{cases}
\end{equation}
$g$ is a biasymptotically quasiperiodic solution associated to $(X, X_{0, \pm}\varphi_{0, \pm})$. 
\end{proposition}
\begin{proof}
Let $\psi_{t_0, X}^t$ be the flow at time $t$ with initial time $t_0$ of $X$. If $\varphi^t_+$ is a positive $C^\sigma$-asymptotic KAM torus associated to $(X, X_{0,+} \varphi_{0,+})$, then
\begin{equation}
\label{congjprop+}
\psi^t_{t_0, X}\circ \varphi_+^{t_0}(q) = \varphi_+^t(q +\omega_+(t - t_0))
\end{equation}
for all $t$, $t_0 \in [0, +\infty)$ and $q \in \T^n$, we refer to~\cite{CdlL15} or~\cite{Sca22a}. Moreover, by~\eqref{Hyp1} and the latter
\begin{equation}
\label{propintro1}
\lim_{t \to +\infty}|\psi^t_{0, X} \circ \varphi_+^0(q) - \varphi_{0,+}(q + \omega_+t)| = 0
\end{equation}
for all $q \in \T^n$. Similarly, 
\begin{align}
\label{congjprop-}
&\psi^t_{t_0, X}\circ \varphi_-^{t_0}(q) = \varphi_-^t(q +\omega_-(t - t_0)), \\
\label{propintro2}
& \lim_{t \to -\infty}|\psi^t_{0, X} \circ \varphi_-^0(q) - \varphi_{0,-}(q + \omega_-t)| = 0.
\end{align}
for all $t$, $t_0 \in (-\infty, 0]$ and $q \in \T^n$. 

Thanks to~\eqref{varphi0+-},~\eqref{congjprop+} and~\eqref{congjprop-}, the curve $g(t)$ defined in~\eqref{g} is equal to 
\begin{equation*}
g(t)=\begin{cases}   \psi^t_{0, X}\circ \varphi_+^{0}(q_+) \quad \mbox{for all $t \ge 0$}\\
\psi^t_{0, X}\circ \varphi_-^{0}(q_-)  \quad \mbox{for all $t \le 0$},
\end{cases}
\end{equation*}
and hence, by~\eqref{propintro1} and~\eqref{propintro2}, it is a biasymptotically quasiperiodic solution associated to $(X, X_{0, \pm}\varphi_{0, \pm})$. 
\end{proof}

\section{Proof of Theorem \ref{ThmInt} assuming Theorem \ref{Thm3BM}}\label{ProofBM1}

In this section, we assume Theorem \ref{Thm3BM} and we deduce Theorem \ref{ThmInt}. First, we introduce the following well-known property. 

\begin{proposition}
\label{propdiAbed}
Given $r >0$ and $0<\epsilon <\delta< r$, let $\phi$ be a map 
\begin{equation*}
\phi : B_r \longrightarrow  B_{r + \delta}
\end{equation*}
of class $C^1$ such that $|\phi - \mathrm{Id}|_{C^1}< \epsilon$. Then, for $\epsilon$ small enough, $\phi$ is a diffeomorphism onto its image and
\begin{equation*}
B_{r-\delta} \subset \phi(B_r).
\end{equation*}
\end{proposition}

For all $p_0 \in B_{3 \over 4}$, we let $p = p_0 + I$ and we expand the Hamiltonian~\eqref{H1BM} around $p_0$ so that 
\begin{eqnarray*}
h(p) &=& h(p_0) + \partial_p h(p_0) \cdot I +  \int_0^1 (1 -\tau) \partial^2_p h(p_\tau)d\tau \cdot I^2\\
f(q,p,t) &=& f(q, p_0, t) + \partial_p f(q, p_0,t) \cdot I +  \int_0^1 (1 -\tau)\partial^2_p f(q, p_\tau, t) d\tau\cdot I^2
\end{eqnarray*}
where $p_\tau = p_0 + \tau I$ and $I \in B_{1 \over 4}$. For any $p_0 \in B_{3 \over 4}$, we define 
\begin{eqnarray*}
e(p_0)  &=& h(p_0)\\
\omega(p_0) &=& \partial_p h(p_0)\\
a(q, t ; p_0) &=& f(q, p_0, t)\\
b(q, t ; p_0)&=& \partial_p f(q,p_0,t)\\
m(q,I,t ; p_0)  &=& \int_0^1 (1-\tau)\left(\partial_p^2 h(p_\tau) + \partial_p^2 f(q, p_\tau, t)\right) d\tau\\
&=& \int_0^1 (1-\tau)\partial_p^2 H (q, p_\tau, t) d\tau,
\end{eqnarray*}
for all $(q,I,t) \in \T^n \times B_{1 \over 4} \times \R$. Writing $\theta$ instead of $q$ for the angular variables, we can rewrite $H$ in the following form as a family of Hamiltonians parametrized by $p_0 \in B_{3 \over 4}$,
\begin{align}
\label{HthmCpar}
&H:\T^n \times  B_{1 \over 4} \times \R \times B_{3 \over 4} \longrightarrow \R \nonumber\\
&H(\theta,I,t;p_0) = e(p_0) + \omega(p_0) \cdot I \\
&\hspace{25mm} + a(\theta, t; p_0) + b(\theta,t; p_0) \cdot I + m(\theta,I,t; p_0) \cdot I^2. \nonumber
\end{align}
In addition, we consider the following family of Hamiltonians
\begin{equation*}
\tilde h(\theta,I,t;p_0) = e(p_0) + \omega(p_0) \cdot I + m(\theta,I,t; p_0) \cdot I^2
\end{equation*}
for all $(\theta,I,t;p_0) \in \T^n \times  B_{1 \over 4} \times \R \times B_{3 \over 4}$. For each fixed $p_0 \in B_{3 \over 4}$, $X_H$ and $X_{\tilde h}$ satisfy~\eqref{hyp1asymKAM}. Furthermore, for all fixed $p_0 \in B_{3 \over 4}$, it is obvious that $\tilde h$ has an invariant torus supporting quasiperiodic dynamics of frequency vector $\omega(p_0)$. 

It is straightforward to verify that the Hamiltonian $H$ defined by~\eqref{HthmCpar} satisfies the hypotheses of Theorem \ref{Thm3BM}. Then, there exists a family of positive $C^\sigma$-asymptotic KAM tori 
\begin{equation*}
\psi^t_+ : \T^n \times B_{3 \over 4} \longrightarrow \T^n \times B_{1 \over 4}
\end{equation*}
associated to $(X_H, X_{\tilde h}, \psi_0)$ and a family of negative $C^\sigma$-asymptotic KAM tori 
\begin{equation*}
\psi^t_- : \T^n \times B_{3 \over 4} \longrightarrow \T^n \times B_{1 \over 4}
\end{equation*}
associated to $(X_H, X_{\tilde h}, \psi_0)$, where $\psi_0$ is the family of trivial embeddings introduced by~\eqref{trivialembeddingBM}. Moreover, we have 
\begin{equation*}
\sup_{t \ge 0} |\psi^t_+ - \psi_0|_{C^1} < C_0 \varepsilon, \quad \sup_{t \le 0} |\psi^t_- - \psi_0|_{C^1} < C_0 \varepsilon,
\end{equation*}
where $C_0$ is a constant depending on $n$, $l$, $\Upsilon$ and $|\omega|_{C^1}$.

Therefore, by the latter, there exist $u^t_\pm$, $v^t_\pm :\T^n \times B_{3 \over 4} \to \R^n$ such that we can rewrite $ \psi^t_+$ and $ \psi^t_-$ in the following form
\begin{equation*}
 \psi^t_\pm(\theta,p_0) =(\theta + u_\pm^t(\theta, p_0), v_\pm^t(\theta, p_0))
\end{equation*}
for all $(\theta, p_0) \in \T^n \times B_{3 \over 4}$ with
\begin{equation*}
\sup_{t \ge 0} |u_+^t|_{C^1} < C_0 \varepsilon, \quad \sup_{t \ge 0 } |v_+^t|_{C^1} < C_0 \varepsilon,\quad \sup_{t \le 0} |u_-^t|_{C^1} < C_0 \varepsilon, \quad \sup_{t \le 0 } |v_-^t|_{C^1} < C_0 \varepsilon.
\end{equation*}

By construction, an orbit $(\theta(t), I(t))$ for the previous Hamiltonian at the parameter value $p_0 \in B_{3 \over 4}$ translates into a trajectory $(q(t), p(t)) = (\theta(t), p_0 + I(t))$ for the Hamiltonian in $(q,p)$-coordinates. Then, letting
\begin{equation}
\label{trivialembeddingproofThmIntBM}
 \varphi_0 : \T^n \times B_{3 \over 4} \longrightarrow \T^n \times B_1, \quad \varphi_0(q,p_0) = (q, p_0),
\end{equation}
the following family of maps
\begin{align*}
&\varphi_\pm^t : \T^n \times B_{3 \over 4} \longrightarrow \T^n \times B_1, \\
&\varphi_\pm^t(q,p_0) = (q + u_\pm^t(q, p_0), p_0 + v_\pm^t(q, p_0))
\end{align*}
is a family of positive (resp. negative) $C^\sigma$-asymptotic KAM tori associated to $(X_H, X_{\tilde h}, \varphi_0)$. In other words, for all $p_0 \in B_{3 \over 4}$, $\varphi^t_{+p_0}$ (resp. $\varphi^t_{-p_0}$) is a positive (resp. negative) $C^\sigma$-asymptotic KAM torus associated to $(X_H, X_{\tilde h}, \varphi_{0,p_0})$.
Thanks to Proposition \ref{propdiAbed},
\begin{equation*}
\T^n \times B_{1 \over 2} \subset \varphi_\pm^0(\T^n \times B_{3 \over 4}).
\end{equation*}
This concludes the proof of the theorem because, for all $(q, p_0) \in \T^n \times B_{1 \over 2}$, there exist $(q_+, p_{0+})$,  $(q_-, p_{0-}) \in \T^n \times B_{3 \over 4}$ such that 
\begin{equation*}
\varphi_+^0(q_+, p_{0+}) = (q,p_0) = \varphi_-^0(q_-, p_{0-}).
\end{equation*}
Then, by Proposition \ref{proplink2BM} there exists a biasymptotically quasiperiodic solution $g(t)$ associated to $(X_H, X_{\tilde h}, \varphi_{0, p_{0\pm}})$ such that $g(0) = (q,p_0)$.

\section{Proof of  Theorem \ref{Thm3BM}}\label{ProofBM2}

The proof of Theorem \ref{Thm3BM} is very similar to what we did in~\cite{Sca22a}. It contains some modifications, especially in the section dedicated to the homological equation. The Hamiltonian in~\eqref{H3BM} consists of a family of Hamiltonians parametrized by $p_0 \in B_{3 \over 4}$. This section aims to prove the existence of a positive $C^\sigma$-asymptotic KAM torus for each $p_0 \in B_{3 \over 4}$. Similarly, we have the claim concerning the existence of a family of negative $C^\sigma$-asymptotic KAM tori. Moreover, we need that these families satisfy~\eqref{stimethm3BM}.

\subsection{Outline of the proof of Theorem \ref{Thm3BM}}
We are looking for a family of positive $C^\sigma$-asymptotic KAM tori $\psi^t$ associated to $(X_H, X_{\tilde h}, \psi_0)$, where we drop the subscript $+$ in order to obtain a more elegant form. Here, $H$ is the Hamiltonian in~\eqref{H3BM} and $\psi_0$ is the following family of trivial embeddings
\begin{equation*}
\psi_0 :\T^n \times B_{3 \over 4} \to \T^n \times B_{1 \over 4}, \quad \psi_0(\theta, p_0) = (\theta, 0). 
\end{equation*}
For the sake of clarity, this means that, for all $p_0 \in B_{3 \over 4}$, we are looking for a positive $C^\sigma$-asymptotic KAM torus $\psi^t_{p_0}$ associated to $(X_H, X_{\tilde h}, \psi_{0, p_0})$.

More specifically, we are searching for $u$, $v : \T^n \times \R^+ \times B_{3 \over 4} \to \R^n$ such that 
\begin{equation*}
\psi(\theta, t;p_0) = (\theta + u(\theta, t;p_0), v(\theta, t;p_0))
\end{equation*}
and,  for all fixed $p_0 \in B_{3 \over 4}$, $\psi$, $u$ and $v$ satisfy the following conditions
\begin{align}
\label{XHOPBM}
&X_H(\psi(\theta, t;p_0), t;p_0) - \partial_\theta \psi(\theta, t;p_0) \omega(p_0) - \partial_t \psi(\theta, t;p_0)=0 \\
\label{XHOPBM2}
& \lim_{t \to +\infty}|u_{p_0}^t|_{C^\sigma}=0, \quad \lim_{t \to +\infty}|v_{p_0}^t|_{C^\sigma}=0, 
\end{align}
for all $(\theta, t) \in \T^n \times \R^+$.

This proof relies on the implicit function theorem. For this reason, we define a suitable functional $\mathcal{F}$ given by~\eqref{XHOPBM}. First, let us introduce the following notation $ \partial_I \Big(m(\theta, I, t; p_0) \cdot I^2 \Big) = \bar m(\theta, I, t; p_0) I $ with 
\begin{equation*}
 \bar m(\theta, I, t; p_0)I = \left(\int_0^1  \partial_p^2 H(\theta ,p_0 + \tau I, t) d\tau \right) I 
\end{equation*}
for all $(\theta, I, t ;p_0) \in \T^n \times B_{1\over 4} \times \R^+ \times B_{3 \over 4}$. Concerning the definition of the functional $\mathcal{F}$, we can see that the Hamiltonian system associated with the Hamiltonian $H$ in~\eqref{H3BM} is equal to 
\begin{equation*}
X_H(\theta,I,t;p_0)  = \begin{pmatrix} \omega(p_0) + b(\theta,t;p_0) + \bar m(\theta,I,t;p_0) I  \\
-\partial_\theta a(\theta,t;p_0)  - \partial_\theta b(\theta,t;p_0)I - \partial_\theta m(\theta,I,t;p_0) I^2\end{pmatrix},
\end{equation*}
for all $(\theta, I, t ;p_0) \in \T^n \times B_{1\over 4} \times \R^+ \times B_{3 \over 4}$. Now, we define
\begin{equation*}
\tilde \psi (\theta, t;p_0) =  (\theta + u(\theta,t; p_0), v(\theta,t; p_0), t ;p_0), \quad \tilde u(\theta,t; p_0) = (\theta + u(\theta,t ;p_0), t ;p_0)
\end{equation*}
for all $(\theta, I, t ;p_0) \in \T^n \times B_{1\over 4} \times \R^+ \times B_{3 \over 4}$. Therefore, the composition of the Hamiltonian system $X_H$ with $\tilde \psi$ can be written as
\begin{equation*}
X_H \circ \tilde \psi = \begin{pmatrix}\omega + b \circ \tilde u + \left(\bar m \circ \tilde \psi\right) v  \\
-\partial_\theta a\circ \tilde u - \left(\partial_\theta b \circ \tilde u \right) v  - \left(\partial_\theta m \circ \tilde \psi\right) \cdot v^2\end{pmatrix}.
\end{equation*}
The latter is composed of sums and products of functions defined on $(\theta, t; p_0) \in \T^n  \times \R^+ \times B_{3 \over 4}$, we have omitted the arguments $(\theta, t; p_0)$ in order to achieve a more elegant form. We keep this notation for the rest of the proof.

On the other side, 
\begin{equation*}
\partial_\theta \psi(\theta, t;p_0) \omega(p_0) + \partial_t \psi(\theta, t;p_0)= \begin{pmatrix} \omega(p_0) + \partial_\theta u(\theta, t;p_0)\omega(p_0) + \partial_t u (\theta, t;p_0) \\
\partial_\theta v(\theta, t;p_0)\omega(p_0) + \partial_t v (\theta, t;p_0)\end{pmatrix}
\end{equation*}
for all $(\theta, t;p_0) \in \T^n \times \R^+ \times B_{3 \over 4}$. Now, we define
\begin{align*}
&\nabla_{\theta t} u (\theta,t; p_0)\Omega(p_0) = \partial_\theta u(\theta,t; p_0) \omega(p_0) + \partial_t u(\theta,t; p_0),\\
&\nabla_{\theta t} v (\theta,t; p_0)\Omega(p_0) = \partial_\theta v(\theta,t; p_0) \omega(p_0) + \partial_t v(\theta,t; p_0).
\end{align*}
for all $(\theta, t;p_0) \in \T^n \times \R^+ \times B_{3 \over 4}$. Then, we can rewrite~\eqref{XHOPBM} in the following form 
\begin{eqnarray}
\label{XHOP2BM}
\begin{pmatrix}b \circ \tilde u + \left(\bar m \circ \tilde \psi  \right) v - \left(\nabla_{\theta t} u \right)\Omega \\
-\partial_\theta a\circ \tilde u  - \left(\partial_\theta b \circ \tilde u \right) v  - \left(\partial_\theta m \circ \tilde \psi  \right) \cdot v^2 - \left(\nabla_{\theta t} v \right)\Omega 
\end{pmatrix} = \begin{pmatrix} 0\\0\end{pmatrix},
\end{eqnarray}
where we have omitted the arguments $(\theta, t; p_0)$. Over suitable Banach spaces, which we will specify later, we consider the following functional 
\begin{equation*}
\mathcal{F}(a, b, m, \bar m,u, v) = (F_1(b, \bar m,u,v), F_2(a,b, m,u,v)),
\end{equation*}
such that
\begin{eqnarray*}
F_1(b, \bar m,u,v) &=& b \circ \tilde u + \left(\bar m \circ \tilde \psi \right) v - \left(\nabla_{\theta t} u \right)\Omega,\\
F_2(a,b, m, u,v) &=& \partial_\theta a\circ \tilde u  + \left(\partial_\theta b \circ \tilde u \right) v  + \left(\partial_\theta m \circ \tilde \psi \right) \cdot v^2 + \left(\nabla_{\theta t} v \right)\Omega .
\end{eqnarray*}
We point out that $\mathcal{F}$ is obtained by~\eqref{XHOP2BM}. Thus, we can reformulate our problem in the following form. For fixed $m$ and $\bar m$ in a suitable Banach space and for $(a,b)$ sufficiently close to $(0,0)$, we are looking for $u$, $v : \T^n \times \R^+ \times B_1 \to \R^n$ satisfying~\eqref{XHOPBM2} such that $\mathcal{F} (a,b,m, \bar m,u,v) = 0$.

Concerning the associated linearized problem, the differential of $\mathcal{F}$ with respect to the variables $(u,v)$ calculated in $(0,0,m, \bar m,0,0)$ is equal to 
\begin{equation*}
D_{(u,v)} \mathcal{F} (0,0,m, \bar m,0,0) (\hat u, \hat v) = (\bar m_0 \hat v - (\nabla_{\theta t} \hat u) \Omega, (\nabla_{\theta t} \hat v) \Omega),
\end{equation*}
where, for all $(\theta, t,p_0) \in \T^n \times B_{1 \over 4} \times B_{3 \over 4}$, we let $\bar m_0(\theta,t ;p_0) = \bar m_0(\theta,0,t ;p_0)$. 

In the following four sections, we prove Theorem \ref{Thm3BM}. First, we introduce suitable Banach spaces on which the functional $\mathcal{F}$ is defined (Section \ref{PSTh2}). In Section \ref{HMTh2}, we solve the homological equation, which is the key to proving that the latter operator is invertible. Then, we verify that $\mathcal{F}$ satisfies the hypotheses of the implicit function theorem (Section \ref{RegTh2}). Finally (Section \ref{FineTh2}), we verify the estimates, and we conclude the proof of the theorem.

\subsection{Preliminary Settings}\label{PSTh2}
We begin this section by introducing some notations, spaces of functions and suitable norms. Given $\sigma \ge 1$, we have the following definition
\begin{definition}
Let $\mathcal{B}^+_{\sigma}$ be the space of functions $f$ defined on $\T^n \times \R^+ \times B_{3 \over 4}$ such that $f$, $\partial_{p_0} f \in C(\T^n \times \R^+ \times B_{3 \over 4})$ and $f_{p_0}^t \in C^\sigma(\T^n)$ for all $(t, p_0) \in \R^+ \times B_{3 \over 4}$.
\end{definition}
For all $f \in \mathcal{B}^+_{\sigma}$ and $l >1$, we define the following norms
\begin{align*}
&|f|^+_{\sigma , l} = \sup_{(t, p_0) \in  \R^+ \times B_{3 \over 4}}|f^t_{p_0}|_{C^\sigma}(1 + t^l) + \sup_{(t, p_0) \in  \R^+ \times B_{3 \over 4}}|\left(\partial_{p_0} f\right)^t_{p_0}|_{C^0}(1 + t^{l-1}),\\
&|f|^+_{\sigma , 0} = \sup_{(t, p_0) \in  \R^+ \times B_{3 \over 4}}|f^t_{p_0}|_{C^\sigma} + \sup_{(t, p_0) \in  \R^+ \times B_{3 \over 4}}|\left(\partial_{p_0} f\right)^t_{p_0}|_{C^0}.
\end{align*}
These norms satisfy the properties in Proposition \ref{normpropertiesBM}. As one can expect, we define the following subset of $\mathcal{B}^+_{\sigma}$

\begin{definition}
Given $\sigma \ge 1$ and an integer $k \ge 0$, we define $\mathcal{\bar B}^+_{\sigma, k}$ the space of functions $f$ such that 
\begin{equation*}
f \in \mathcal{B}^+_{\sigma +k}, \hspace{2mm} \mbox{and} \hspace{2mm} \partial^i_q f \in \mathcal{B}^+_{\sigma + k -i}
\end{equation*}
for all $0 \le i \le k$. 
\end{definition} 
We conclude this part of settings with the following norm. For all $f \in \mathcal{\bar B}^+_{\sigma, k}$ and $l>1$, we define
\begin{equation*}
\left \|f \right \|^+_{\sigma, k, l} = \max_{0 \le i \le k}|\partial_q^i f|^+_{\sigma+k-i, l}.
\end{equation*} 

Now, let $\sigma \ge 1$ and $l >1$ be the positive parameters in~\eqref{H3BM}. We consider the following Banach spaces $\left(\mathcal{A}, |\cdot |\right)$, $\left(\mathcal{B}, |\cdot |\right)$, $\left(\mathcal{U}, |\cdot |\right)$, $\left(\mathcal{V}, |\cdot |\right)$, $\left(\mathcal{Z}, |\cdot |\right)$ and $\left(\mathcal{G}, |\cdot |\right)$ 
\vspace{5mm}
\begin{eqnarray*}
\mathcal{A} &=& \Big\{a : \T^n \times \R^+ \times B_{3 \over 4} \to \R  \hspace{1mm}| \hspace{1mm} a \in \mathcal{\bar B}^+_{\sigma, 2} \hspace{1mm} \mbox{and} \hspace{1mm} |a| = |a|^+_{\sigma+2,0} + \left \|\partial_\theta a\right\|^+_{\sigma,1, l+2} < \infty\Big\}\\
\mathcal{B} &=& \Big\{b : \T^n \times \R^+ \times B_{3 \over 4}  \to \R^n  \hspace{1mm}| \hspace{1mm} b \in \mathcal{\bar B}^+_{\sigma, 2,}, \hspace{1mm} \mbox{and} \hspace{1mm} |b| = \left \|b\right\|^+_{\sigma,2, l+1} < \infty\Big\}\\
\mathcal{U} &=& \Big\{u :\T^n \times \R^+ \times B_{3 \over 4}  \to \R^n  \hspace{1mm}| \hspace{1mm} u ,  \left( \nabla_{\theta t} u\right) \Omega \in \mathcal{B}^+_{\sigma}  \\
&& \mbox{and} \hspace{1mm} |u| =\max \{|u|^+_{\sigma,l},| \left( \nabla_{\theta t} u \right) \Omega|^+_{\sigma,l+1} \} < \infty\Big\}\\
\mathcal{V} &=& \Big\{v : \T^n \times \R^+ \times B_{3 \over 4}  \to \R^n  \hspace{1mm}| \hspace{1mm} v,  \left( \nabla_{\theta t} v\right) \Omega \in \mathcal{B}^+_{\sigma}  \\
&& \mbox{and} \hspace{1mm} |v| =\max \{|v|^+_{\sigma,l+1},| \left( \nabla_{\theta t} v \right) \Omega|^+_{\sigma,l+2} \} < \infty\Big\}\\
\mathcal{Z} &=& \Big\{z : \T^n \times \R^+ \times B_{3 \over 4}  \to \R^n  \hspace{1mm}| \hspace{1mm} z \in \mathcal{B}^+_{\sigma},  \hspace{1mm} \mbox{and} \hspace{1mm} |z| =|z|^+_{\sigma,l+1} < \infty\Big\}\\
\mathcal{G} &=& \Big\{g : \T^n \times \R^+ \times B_{3 \over 4}  \to \R^n  \hspace{1mm}| \hspace{1mm} g \in \mathcal{B}^+_{\sigma},  \hspace{1mm} \mbox{and} \hspace{1mm} |g| =|g|^+_{\sigma,l+2} < \infty\Big\}
\end{eqnarray*} 

Similarly to what we did in~\cite{Sca22a}, verifying that the previous normed spaces are Banach spaces is straightforward. Let $M_n$ be the set of the $n$-dimensional matrices and $\Upsilon \ge 1$ the positive parameter in~\eqref{H3BM}. We introduce another Banach space $\left(\mathcal{M}, |\cdot |\right)$, such that 
\begin{eqnarray*}
\mathcal{M} &=& \Big\{m : \T^n \times B_{1 \over 4} \times \R^+ \times B_{3 \over 4} \to M_n  \hspace{1mm}| \hspace{1mm}  \partial^i_{\theta I p_0}m \in C(\T^n \times B_{1 \over 4} \times \R^+ \times B_{3 \over 4}) \hspace{1mm}\\
&&\mbox{for all $0 \le i \le 3$,}  \hspace {2mm}m^t \in C^{\sigma+2}( \T^n \times B_{1 \over 4} \times B_{3 \over 4}) \hspace{1mm} \mbox{for all fixed $t \in \R^+$} \\
&&\mbox{and}\hspace{1mm} |m| = \hspace{1mm}\sup_{t \in \R^+}|m^t|_{C^{\sigma + 2}} \le \Upsilon\Big\}
\end{eqnarray*}

Now, we can define the functional $\mathcal{F}$ introduced in the previous section more precisely.
Let $\mathcal{F}$ be the following functional 
\begin{equation}
\label{F1}
\mathcal{F} :  \mathcal{A} \times \mathcal{B} \times \mathcal{M} \times \mathcal{M} \times  \mathcal{U} \times  \mathcal{V} \longrightarrow \mathcal{Z} \times \mathcal{G}\\
\end{equation}
\begin{equation*}
\mathcal{F}(a, b,m , \bar m,  u, v) = (F_1(b, \bar m, u,v), F_2(a,b,m,u,v))
\end{equation*}
with 
\begin{eqnarray*}
F_1(b,\bar m ,u,v) &=& b \circ \tilde u + \left(\bar m \circ \tilde \psi \right) v - \left(\nabla_{\theta t} u \right)\Omega,\\
F_2(a,b, m,u,v) &=& \partial_\theta a\circ \tilde u  + \left(\partial_\theta b \circ \tilde u \right) v  + \left(\partial_\theta m \circ \tilde \psi \right) \cdot v^2 + \left(\nabla_{\theta t} v \right)\Omega .
\end{eqnarray*}

\subsection{Homological equation}\label{HMTh2}

Before analyzing the homological equation, let us prove the following estimates
\begin{lemma}
\label{lemmastimeBMjo}
Given $m >1$,
\begin{equation}
\label{stimascemaBMsul}
\int_t^{+\infty}{1 \over 1 + \tau^m}d\tau \le {C(m)\over  1 + t^{m-1}}, \quad \int_t^{+\infty}{\tau -t \over 1 + \tau^{m+1}}d\tau \le {C(m)\over  1 + t^{m-1}}
\end{equation}
for all $t \ge 0$ and some constants $C(m)$ depending on $m$. 
\end{lemma}
\begin{proof}
We define the following function $f_m :\R^+ \longrightarrow \R$ such that
\begin{equation*}
f_m(t) = \left(1 + t^{m-1}\right)\int_t^{+\infty}{1 \over 1 + \tau^m}d\tau.
\end{equation*}
It is straightforward to verify that $f$ is continuous. We will prove the existence of a constant $C(m)$, depending on $m$, such that $f_m(t) \le C(m)$ for all $t \ge 0$, which implies the first estimate in~\eqref{stimascemaBMsul}. It suffices to prove that there exists $\lim_{t \to +\infty} f_m(t)$, and it is finite. Thanks to l'Hôpital's rule
\begin{eqnarray*}
\lim_{t \to +\infty} f_m(t) = \lim_{t \to +\infty} {{d \over dt}\int_t^{+\infty}{1 \over 1 + \tau^m}d\tau \over {d \over dt} {1 \over 1 + t^{m-1}}} = \lim_{t \to +\infty} {\left(1 + t^{m-1}\right)^2 \over (m-1)t^{m-2}\left(1 + t^m\right)} = {1 \over m-1}.
\end{eqnarray*}
Concerning the second inequality in~\eqref{stimascemaBMsul}, similarly to the previous case, we define the following function $g_m :\R^+ \longrightarrow \R$ such that 
\begin{equation*}
g_m(t) = \left(1 + t^{m-1}\right)\int_t^{+\infty}{\tau -t \over 1 + \tau^{m+1}}d\tau.
\end{equation*}
One can see that $g_m$ is continuous. We have to verify that there exists $\lim_{t \to +\infty} g_m(t)$, and it is finite. Applying l'Hôpital's rule twice
\begin{eqnarray*}
\lim_{t \to +\infty} g_m(t) &=& \lim_{t \to +\infty} {{d \over dt}\int_t^{+\infty}{\tau -t \over 1 + \tau^{m+1}}d\tau \over {d \over dt} {1 \over 1 + t^{m-1}}} = \lim_{t \to +\infty} {\int_t^{+\infty}{1 \over 1 + \tau^m}d\tau \over {(m-1)t^{m-2} \over \left(1 + t^{m-1}\right)^2}} = \lim_{t \to +\infty} {{d \over dt}\int_t^{+\infty}{1 \over 1 + \tau^m}d\tau \over {d \over dt}{(m-1)t^{m-2} \over \left(1 + t^{m-1}\right)^2}}\\
&=&\lim_{t \to +\infty}{\left(1 + t^{m-1}\right)^4 \over \left(1 + t^{m+1}\right)t^{3m-5}h_m(t)}
\end{eqnarray*}
where 
\begin{equation*}
h_m(t) = (m-1)\left(2 (m-1) \left({1 \over t^{m-1}}+1\right) - (m-2)\left({1 \over t^{m-1}}+1\right)^2\right). 
\end{equation*}
Then, by the latter
\begin{eqnarray*}
\lim_{t \to +\infty} g_m(t) &=&\lim_{t \to +\infty}{\left(1 + t^{m-1}\right)^4 \over \left(1 + t^{m+1}\right)t^{3m-5}h_m(t)} = \lim_{t \to +\infty}{\left({1 \over t^{m-1}} + 1\right)^4 \over \left({1 \over  t^{m+1}} +1\right)h_m(t)} = {1 \over m(m-1)}.
\end{eqnarray*}
\end{proof}

Given $\sigma \ge 1$ and $l>1$, we consider the following equation in the unknown $\varkappa : \T^n \times \R^+ \times B_{3 \over 4} \to \R$
\begin{equation}
\label{HEBM}
\begin{cases}
 \omega(p_0) \cdot \partial_q \varkappa(\theta, t ; p_0) + \partial_t \varkappa(\theta, t ; p_0) = g(\theta, t ; p_0),\\
g \in \mathcal{B}^+_{\sigma}, \quad |g|^+_{\sigma, l + 1} < \infty,\\
\omega : B_{3 \over 4} \longrightarrow \R^n, \quad \omega \in C^1(B_{3 \over 4}).\\
\end{cases}
\tag{$HE_A$}
\end{equation}

\begin{lemma}[\textbf{Homological Equation}]
\label{homoeqlemmaBM} 
There exists a unique solution $\varkappa \in \mathcal{B}^+_{\sigma}$ of~\eqref{HEBM} such that, for all fixed $p_0 \in B_{3 \over 4}$,
\begin{equation}
\label{HEasycond}
\lim_{t \to +\infty} |\varkappa^t_{p_0}|_{C^0} = 0.
\end{equation}
Moreover, 
\begin{equation}
\label{varkappaestimateBM1}
|\varkappa|^+_{\sigma, l} \le C(l, |\omega|_{C^1}) |g|^+_{\sigma, l+1}
\end{equation}
for a suitable constant $C(l, |\omega|_{C^1})$ depending on $n$, $l$ and $|\omega|_{C^1}$. 
\end{lemma}

\begin{proof}
\textit{Existence}: This first part is extremely similar to what we did in~\cite{Sca22a}. For this reason, we refer to~\cite{Sca22a} for more detailed proof.

  Let us define the following transformation 
\begin{equation*}
h: \T^n \times \R^+ \times B_{3 \over 4}  \to  \T^n \times \R^+ \times B_{3 \over 4}, \quad h(\theta, t ; p_0) = (\theta -  \omega(p_0) t, t; p_0).
\end{equation*}
We claim that it is enough to prove the first part of this lemma for the much simpler equation
\begin{equation*}
\partial_t \kappa(\theta, t ; p_0) = g(\theta + \omega(p_0) t, t ; p_0).
\end{equation*}
This is because if $\kappa$ is a solution of the above equation satisfying~\eqref{HEasycond}, then $\chi = \kappa \circ h$ is a solution of (\ref{HEBM}) satisfying~\eqref{HEasycond} and viceversa. 

The unique solution of the above equation satisfying the asymptotic condition is
\begin{equation*}
\kappa(\theta, t ; p_0) = -\int_t^{+\infty} g(\theta + \omega(p_0) \tau, \tau ; p_0)d\tau
\end{equation*}
and hence composing $k$ with $h$
\begin{equation}
\label{chiBM1}
\varkappa(\theta, t ; p_0) = \kappa \circ h(\theta, t ; p_0) =-\int_t^{+\infty} g(\theta+ \omega(p_0)(\tau - t),\tau ,p_0)d\tau
\end{equation}
is the unique solution of (\ref{HEBM}) that we are looking for. 

\textit{Regularity and Estimates}: $g \in \mathcal{B}^+_{\sigma}$ implies $\kappa \in \mathcal{B}^+_{\sigma}$ and thus $\varkappa = \kappa \circ h \in \mathcal{B}^+_{\sigma}$. Now, we have to verify the estimate~\eqref{varkappaestimateBM1}. By~\eqref{chiBM1} and Lemma \ref{lemmastimeBMjo}
\begin{equation*}
|\varkappa_{p_0}^t|_{C^\sigma} \le \int_t^{+\infty}|g^\tau_{p_0}|_{C^\sigma} d\tau \le |g|^+_{\sigma,l+1}\int_t^{+\infty}{1 \over 1 + \tau^{l+1}}d\tau \le C(l){|g|^+_{\sigma,l+1} \over 1 + t^l},
\end{equation*}
for all fixed $(t, p_0) \in  \R^+ \times B_{3 \over 4}$. Multiplying both sides of the latter by $1 + t^l$ and taking the sup for all $ \R^+ \times B_{3 \over 4}$, we obtain
\begin{equation*}
\sup_{(t, p_0) \in  \R^+ \times B_{3 \over 4}}|\varkappa_{p_0}^t|_{C^\sigma}(1 + t^l) \le C(l)|g|^+_{\sigma,l+1}.
\end{equation*}
It remains to estimate the second member of the norm $|\cdot|^+_{\sigma, l}$.  The partial derivate of $\varkappa$ with respect to $p_0$ is equal to
\begin{eqnarray*}
\partial_{p_0} \varkappa(\theta, t ; p_0) &=&-\int_t^{+\infty}  \partial_{p_0} \omega(p_0)\partial_\theta g(\theta+ \omega(p_0)(\tau - t),\tau ; p_0) (\tau - t) d\tau \\
&-&\int_t^{+\infty}  \partial_{p_0} g(\theta+ \omega(p_0)(\tau - t),\tau ; p_0)d\tau.
\end{eqnarray*}
Then, thanks to Lemma \ref{lemmastimeBMjo}, we can estimate the norm $C^0$ on the left-hand side of the latter as follows
\begin{eqnarray*}
|\partial_{p_0} \varkappa_{p_0}^t|_{C^0} &\le& \int_t^{+\infty}|g^\tau_{p_0}|_{C^1}|\omega|_{C^1}(\tau -t) + |\partial_{p_0} g^\tau|_{C^0} d\tau\\
&\le&|g|^+_{1, l+1} |\omega|_{C^1}\int_t^{+\infty}{(\tau -t) \over 1+\tau^{l+1}}d\tau + |g|^+_{1, l+1}\int_t^{+\infty}{1 \over 1+\tau^l}d\tau,\\
&\le& C(l, |\omega|_{C^1}){|g|^+_{1, l+1} \over 1 + t^{l-1}}.
\end{eqnarray*}
Similarly to the previous case, by multiplying both sides of the previous inequality by $1 + t^{l-1}$ and taking the sup for all $ \R^+ \times B_{3 \over 4}$, we have
\begin{equation*}
\sup_{(t, p_0) \in  \R^+ \times B_1}|\partial_{p_0}\varkappa_{p_0}^t|_{C^0}(1 + t^{l-1}) \le C(l, |\omega|_{C^1})|g|^+_{1, l+1}.
\end{equation*}
Moreover, reminding that $|g|^+_{1, l+1} \le |g|^+_{\sigma, l+1}$, we conclude the proof of this lemma. 
\end{proof}

\subsection{Regularity of $\mathcal{F}$}\label{RegTh2}

This section is dedicated to the functional $\mathcal{F}$ defined by~\eqref{F1}. We want to prove that $\mathcal{F}$ satisfies the hypothesis of the implicit function theorem. By Proposition \ref{normpropertiesBM}, $\mathcal{F}$ is well defined, continuous, differentiable with respect to the variables $(u,v)$ and this differential $D_{(u,v)}\mathcal{F}$ is continuous. As we have already seen, $D_{(u,v)}\mathcal{F}$ calculated in $(0,0,m,\bar m,0,0)$ is equal to
\begin{equation}
\label{DFBM}
D_{(u,v)}  \mathcal{F}(0,0,m,\bar m,0,0)(\hat u, \hat v)  = (\bar m_0 \hat v - (\nabla_{\theta t} \hat u) \Omega, (\nabla_{\theta t} \hat v) \Omega).
\end{equation}
It remains to verify that the latter is invertible for all fixed $m$, $\bar m \in\mathcal{M}$. 
\begin{lemma}
\label{LemmaInvBM}
For all $(z,g) \in \mathcal{Z} \times \mathcal{G}$ there exists a unique $(\hat u, \hat v) \in \mathcal{U} \times \mathcal{V}$ such that 
\begin{equation*}
D_{(u,v)} \mathcal{F} (0,0,m, \bar m, 0,0) (\hat u, \hat v) = (z,g).
\end{equation*}
Moreover, for a suitable constant $\bar C$ depending on $n$, $l$ and $|\omega|_{C^1}$
\begin{equation}
\label{uvLemmaInvBM}
|\hat u| \le \bar C \left(\left|\bar m_0\right|^+_{\sigma, 0} |g|^+_{\sigma,l+2}+ \left|z \right|^+_{\sigma, l+1}\right), \quad |\hat v|  \le \bar C |g|^+_{\sigma, l+1},
\end{equation}
where we recall that $|u| =\max \{|u|^+_{\sigma,l},| \left( \nabla_{\theta t} u \right) \Omega|^+_{\sigma,l+1}\}$ \newline and $ |v| =\max \{|v|^+_{\sigma,l+1},| \left( \nabla_{\theta t} v \right) \Omega|^+_{\sigma,l+2} \} $.
\end{lemma}

\begin{proof}
The key of the proof is Lemma \ref{homoeqlemmaBM}. By~\eqref{DFBM}, the proof consists in searching for the unique solution to the following system
\begin{equation}
\label{LPBM2}
\begin{cases}
\bar m_0 \hat v - (\nabla_{\theta t} \hat u) \Omega = z \\
(\nabla_{\theta t} \hat v) \Omega = g.
\end{cases}
\end{equation}
Thanks to Lemma \ref{homoeqlemmaBM}, a unique solution $\hat v$ for the last equation of the above system exists and
\begin{equation}
\label{hatvBM}
|\hat v|^+_{\sigma,l+1} \le C(|\omega|_{C^1}, l) |g|^+_{\sigma,l+2}.
\end{equation}
Moreover, by $\left|\left( \nabla_{\theta t} \hat v \right) \Omega\right|^+_{\sigma,l+2} = |g|^+_{\sigma,l+2}$, we obtain the following estimate
\begin{equation*}
|\hat v|  =\max \{|\hat v|^+_{\sigma,l+1},| \left( \nabla_{\theta t} \hat v \right) \Omega|^+_{\sigma,l+2} \}  \le C(|\omega|_{C^1}, l)|g|^+_{\sigma, l+1}.
\end{equation*}
Now, it remains to solve the first equation of the system~\eqref{hatvBM} where $\hat v$ is known. We can rewrite this equation in the following form
\begin{equation}
\label{LPBM2}
 (\nabla_{\theta t} \hat u) \Omega = \bar m_0 \hat v - z.
\end{equation}
By Proposition \ref{normpropertiesBM} and~\eqref{hatvBM},  we can estimate the norm $|\cdot|_{\sigma, l+1}$ on the right-hand side of the latter as follows
\begin{equation*}
\left|\bar m_0 \hat v - z \right|^+_{\sigma, l+1} \le \left|\bar m_0\right|^+_{\sigma, 0} \left|\hat v\right|^+_{\sigma, l+1} + \left|z \right|^+_{\sigma, l+1} \le C(|\omega|_{C^1}, l) \left|\bar m_0\right|^+_{\sigma, 0} |g|^+_{\sigma,l+2}+ \left|z \right|^+_{\sigma, l+1}.
\end{equation*}
This implies 
\begin{equation*}
\left|\left( \nabla_{\theta t} \hat u \right) \Omega\right|^+_{\sigma,l+1} = \left|\bar m_0 \hat v - z \right|^+_{\sigma, l+1} \le C(|\omega|_{C^1}, l) \left|\bar m_0\right|^+_{\sigma, 0} |g|^+_{\sigma,l+2}+ \left|z \right|^+_{\sigma, l+1}
\end{equation*}
and thanks to Lemma \ref{homoeqlemmaBM}, a unique solution of~\eqref{LPBM2} exists satisfying
\begin{equation*}
\left|\hat u \right|^+_{\sigma,l}  \le C(|\omega|_{C^1}, l) \left(\left|\bar m_0\right|^+_{\sigma, 0} |g|^+_{\sigma,l+2}+ \left|z \right|^+_{\sigma, l+1}\right).
\end{equation*}
This concludes the proof of this lemma because 
\begin{equation*}
|\hat u| =\max \{|\hat u|^+_{\sigma,l},| \left( \nabla_{\theta t} \hat u \right) \Omega|^+_{\sigma,l+1}\} \le C(|\omega|_{C^1}, l) \left(\left|\bar m_0\right|^+_{\sigma, 0} |g|^+_{\sigma,l+2}+ \left|z \right|^+_{\sigma, l+1}\right).
\end{equation*}
\end{proof}

\subsection{Families of $C^\sigma$-asymptotic tori}\label{FineTh2}

The functional $\mathcal{F}$ satisfies the hypotheses of the implicit function theorem. Then, for $\varepsilon$ small enough, there exists a family of positive $C^\sigma$-asymptotic KAM tori 
\begin{equation*}
\psi^t_+ : \T^n \times B_{3 \over 4} \longrightarrow \T^n \times B_{1 \over 4}
\end{equation*}
associated to $(X_H, X_{\tilde h}, \psi_0)$, where $\psi_0$ is the following family of trivial embeddings
\begin{equation*}
\psi_0 :\T^n \times B_{3\over 4} \to \T^n \times B_{1 \over 4}, \quad \psi_0(\theta, p_0) = (\theta, 0). 
\end{equation*}
It remains to verify the estimates~\eqref{stimethm3BM}. To this end, let us state a quantitative version of the implicit function theorem. 

\begin{theorems}[Implicit Function Theorem]
\label{IFT}
Let $\left(\mathcal{X}, |\cdot |\right)$, $\left(\mathcal{Y}, |\cdot |\right)$ and $\left(\mathcal{Z}, |\cdot |\right)$ be Banach spaces. For some $(x_0, y_0) \in \mathcal{X} \times \mathcal{Y}$ and $\varepsilon$, $\mu >0$, we introduce the following spaces
\begin{equation*}
\mathcal{X}_0 = \{x \in \mathcal{X} : |x- x_0| \le \varepsilon\}, \quad \mathcal{Y}_0 = \{y \in \mathcal{Y} : |y- y_0| \le \mu\}. 
\end{equation*}
We assume that 
\begin{equation*}
\mathcal{F} : \mathcal{X}_0 \times \mathcal{Y}_0 \longrightarrow \mathcal{Z}
\end{equation*}
is continuous and has the property that $D_y\mathcal{F}$ exists and is continuous at each point of $\mathcal{X}_0 \times \mathcal{Y}_0$. Moreover, $D_y\mathcal{F}(x_0, y_0)$ is invertible and 
\begin{align}
\label{IFTQBM1}
&\sup_{x \in \mathcal{X}_0}\left| D_y \mathcal{F}(x_0, y_0)^{-1} \mathcal{F}(x, y_0) \right| \le {\mu \over 2} \\
\label{IFTQBM2}
&\sup_{(x, y) \in \mathcal{X}_0 \times \mathcal{Y}_0}\left \| \mathrm{Id} - D_y \mathcal{F}(x_0, y_0)^{-1} D_y\mathcal{F}(x, y) \right \| \le {1 \over 2}
\end{align}
where $\mathrm{Id} \in M_n$ is the identity matrix and $\left \| \cdot \right \|$ stands for the operator norm. Then there exists a unique $g \in C(\mathcal{X}_0, \mathcal{Y}_0)$ such that 
\begin{equation*}
g(x_0) = y_0 \quad \mbox{and} \quad \mathcal{F}(g(x), x)=0
\end{equation*}
for all $x \in \mathcal{X}_0$.
\end{theorems}
\begin{proof}
We refer to~\cite{ChiAM2}.
\end{proof}

We observe that condition~\eqref{IFTQBM1} tells us how to choose $\mu$ as a function of $\varepsilon$. Then~\eqref{IFTQBM2} determines a threshold for $\varepsilon$. In order to conclude the proof of Theorem \ref{Thm3BM}, we have to establish the relation between $\varepsilon$ and $\mu$. 

\begin{lemma}
The estimates~\eqref{stimethm3BM} are satisfied. 
\end{lemma}
\begin{proof}
The proof relies on Lemma \ref{LemmaInvBM} and~\eqref{IFTQBM1}.  For all fixed $\bar m$, $m \in \mathcal{M}$ and for all $(a, b) \in \mathcal{A} \times \mathcal{B}$,
\begin{equation*}
D_{(u,v)} \mathcal{F}(0,0, m, \bar m,0,0)^{-1}\mathcal{F}(a, b, m,\bar m,0,0) = D_{(u,v)} \mathcal{F}(0,0, m, \bar m,0,0)^{-1}\begin{pmatrix} b \\ \partial_\theta a \end{pmatrix}.
\end{equation*}
We want to estimate the right-hand side of the latter. We observe that we can reformulate this problem in terms of estimating the unique solution $(\hat u, \hat v) \in \mathcal{Y}$ of the following system
\begin{equation*}
D_{(u,v)}\mathcal{F}(0,0, m, \bar m,0,0) (\hat u, \hat v) = \begin{pmatrix} b \\ \partial_\theta a \end{pmatrix}.
\end{equation*}
By Lemma \ref{LemmaInvBM}, this solution exists and 
\begin{eqnarray*}
|\hat u| &\le& \bar C \left(\Upsilon |\partial_\theta a|^+_{\sigma, l+2}  + |b|^+_{\sigma, l+1}\right) \le \bar C \Upsilon \varepsilon\\
|\hat v| &\le& \bar C|\partial_\theta a|^+_{\sigma, l+2} \le \bar C \varepsilon,
\end{eqnarray*}
where $\bar C$ is a constant depending on $n$, $l$ and $|\omega|_{C^1}$.
Using the notation of Theorem \ref{IFT}, by~\eqref{uvLemmaInvBM}, we can choose $\mu = 2\bar C \Upsilon \varepsilon$. 

Now, we observe that for all fixed $\bar m$, $m \in \mathcal{M}$
\begin{equation*}
\mathrm{Id} - D_{(u,v)}\mathcal{F}(0,0,m, \bar m, 0,0)^{-1} D_{(u,v)}\mathcal{F}(a,b,m,\bar m,u,v)
\end{equation*}
is continuous with respect to $(a,b,u,v) \in \mathcal{A} \times \mathcal{B}\times \mathcal{U}\times \mathcal{V}$. Then, thanks to Lemma \ref{LemmaInvBM}, one can see that there exists $\varepsilon_0$ depending on $n$, $l$, $\Upsilon$ and $|\omega|_{C^1}$ such that, for all $\varepsilon \le \varepsilon_0$,~\eqref{IFTQBM2} is satisfied.
\end{proof}

We proved the existence of a family of positive $C^\sigma$-asymptotic KAM tori associated to $(X_H, X_{\tilde h}, \psi_0)$ verifying~\eqref{stimethm3BM}. Similarly, we have the existence of a family of negative $C^\sigma$-asymptotic KAM tori associated to $(X_H, X_{\tilde h}, \psi_0)$ verifying~\eqref{stimethm3BM}. This concludes the proof of Theorem \ref{Thm3BM}.

\section{Proof of Theorem \ref{ThmNearInt} assuming Theorem \ref{Thm4BM}}\label{ProofBM3}

Here, we assume Theorem \ref{Thm4BM}, and we prove Theorem \ref{ThmNearInt}. The following well-known property is the Lipschitz version of Proposition \ref{propdiAbed}
\begin{proposition}
\label{propdiAbedLip}
Given $r >0$ and $0<\epsilon <\delta<r$, let $\phi$ be a Lipschitz map 
\begin{equation*}
\phi : B_r \longrightarrow  B_{r + \delta}
\end{equation*}
 such that $|\phi - \mathrm{Id}|_{L(B_r)}< \epsilon$. Then, for $\epsilon$ small enough, $\phi$ is a lipeomorphism onto its image and
\begin{equation*}
B_{r-\delta} \subset \phi(B_{r+\delta}).
\end{equation*}
\end{proposition}
We define $\delta=2C_0\varepsilon$, where $C_0$ is the constant introduced in Theorem \ref{Thm4BM}. Therefore, we consider
\begin{equation*}
D'=B_{1-\delta} \cap D.
\end{equation*}
Now, for all $p_0 \in D'$, we let $p = p_0 + I$. Similarly to the proof of Theorem \ref{ThmInt}, we expand the Hamiltonian $H$ in~\eqref{H2BM} around $p_0$ in such a way that
\begin{eqnarray*}
h(p) &=& h(p_0) + \partial_p h(p_0) \cdot I +  \int_0^1 (1 -\tau) \partial^2_p h(p_\tau)d\tau \cdot I^2\\
R(q,p) &=& R(q, p_0) + \partial_p R(q, p_0) \cdot I +  \int_0^1 (1 -\tau)\partial^2_p R(q, p_\tau, t) d\tau\cdot I^2\\
&=& \int_0^1 (1 -\tau)\partial^2_p R(q, p_\tau, t) d\tau\cdot I^2\\
f(q,p,t) &=& f(q, p_0, t) + \partial_p f(q, p_0,t) \cdot I +  \int_0^1 (1 -\tau)\partial^2_p f(q, p_\tau, t) d\tau\cdot I^2.
\end{eqnarray*}
where $p_\tau = p_0 + \tau I$ and $I \in B_\delta$. For all $p_0 \in D'$, we define
\begin{eqnarray*}
e(p_0)  &=& h(p_0)\\
\omega(p_0) &=& \partial_p h(p_0)\\
a(q, t ; p_0) &=& f(q, p_0, t)\\
b(q, t ; p_0)&=& \partial_p f(q,p_0,t)\\
m(q,I,t ; p_0)  &=& \int_0^1 (1-\tau)\left(\partial_p^2 h(p_\tau) + \partial^2_p R(q, p_\tau, t) + \partial_p^2 f(q, p_\tau, t)\right) d\tau\\
&=& \int_0^1 (1-\tau)\partial_p^2 H (q, p_\tau, t) d\tau
\end{eqnarray*}
for all $(q,I,t) \in \T^n \times B_\delta \times \R$. Writing $\theta$ instead of $q$ for the angular variables, we rewrite the Hamiltonian  $H$, restricted to  $\T^n \times  B_{\delta} \times \R \times D'$, in the following form as a family of Hamiltonians parametrized by $p_0 \in D'$,
\begin{equation}
\label{HThmDproofNotC}
H:\T^n \times  B_{\delta} \times \R \times D' \longrightarrow \R
\end{equation}
such that 
\begin{eqnarray*}
H(\theta,I,t;p_0) &=& e(p_0) + \omega(p_0) \cdot I + a(\theta, t; p_0) + b(\theta,t; p_0) \cdot I + m(\theta,I,t; p_0) \cdot I^2.
\end{eqnarray*}
For all $(\theta,I,t;p_0) \in \T^n \times  B_{\delta} \times \R \times D'$, let $\tilde h$ be the following family of Hamiltonians
\begin{equation*}
\tilde h(\theta,I,t;p_0) = e(p_0) + \omega(p_0) \cdot I + m(\theta,I,t; p_0) \cdot I^2.
\end{equation*}
Obviously, for each fixed $p_0 \in D'$, $X_H$ and $X_{\tilde h}$ satisfy~\eqref{hyp1asymKAM} and, for all fixed $p_0 \in D'$, $\tilde h$ has an invariant torus supporting quasiperiodic dynamics of frequency vector $\omega(p_0)$. 

The Hamiltonian $H$ in~\eqref{HThmDproofNotC} verifies the hypotheses of Theorem \ref{Thm4BM}. Then, similarly to the proof of Theorem \ref{ThmInt}, there exist $u^t_\pm$, $v^t_\pm :\T^n \times D' \to \R^n$ such that, in the $(q,p)$-coordinates, the following family of embeddings
\begin{align*}
&\varphi_+^t : \T^n \times D' \longrightarrow \T^n \times B_1, \\
&\varphi_+^t(q,p_0) = (q + u_+^t(q, p_0), p_0 + v_+^t(q, p_0))
\end{align*}
is a family of positive $C^\sigma$-asymptotic KAM tori associated to $(X_H, X_{\tilde h}, \varphi_0)$ and 
\begin{align*}
&\varphi_-^t : \T^n \times D' \longrightarrow \T^n \times B_1, \\
&\varphi_-^t(q,p_0) = (q + u_-^t(q, p_0), p_0 + v_-^t(q, p_0))
\end{align*}
is a family of negative $C^\sigma$-asymptotic KAM tori associated to $(X_H, X_{\tilde h}, \varphi_0)$, where $\varphi_0$ is the family of trivial embeddings defined by
 \begin{equation*}
 \varphi_0 : \T^n \times D' \longrightarrow \T^n \times B_1, \quad \varphi_0(q,p_0) = (q, p_0).
\end{equation*}
Moreover,
\begin{eqnarray*}
&& \sup_{t \ge 0} |u_+^t|_{L(\T^n \times D')} < C_0 \varepsilon, \hspace{15mm} \sup_{t \ge 0 } |v_+^t|_{L(\T^n \times D')} < C_0 \varepsilon,\\
&& \sup_{t \le 0} |u_-^t|_{L(\T^n \times D')} < C_0 \varepsilon, \hspace{15mm} \sup_{t \le 0 } |v_-^t|_{L(\T^n \times D')} < C_0 \varepsilon,
\end{eqnarray*}
where $C_0$ is a constant depending on $n$, $l$, $\Upsilon$ and $|\omega|_{L(D')}$.

Now, there exist $\tilde u^t_\pm$, $\tilde v^t_\pm :\T^n \times B_{1-\delta} \to \R^n$ such that $\tilde u^t_\pm$, $\tilde v^t_\pm$ extend  $u^t_\pm$, $ v^t_\pm$ without affecting their Lipschitz constant. This means that, 
\begin{align*}
&\hspace{35mm} \tilde u^t_\pm \big|_{\T^n \times D'} =  u^t_\pm, \quad \tilde v^t_\pm \big|_{\T^n \times D'} =  v^t_\pm,\\
&\sup_{t \ge 0} | \tilde u_+^t|_{L(\T^n \times B_{1-\delta})} = \sup_{t \ge 0} |u_+^t|_{L(\T^n \times D')} \quad \sup_{t \ge 0} | \tilde v_+^t|_{L(\T^n \times B_{1-\delta})} = \sup_{t \ge 0} |u_+^t|_{L(\T^n \times D')}. 
\end{align*}
Therefore, letting 
\begin{align*}
&\tilde \varphi_\pm^t : \T^n \times B_{1-\delta} \longrightarrow \T^n \times B_1, \\
&\tilde \varphi_\pm^t(q,p_0) = (q + \tilde u_\pm^t(q, p_0), p_0 + \tilde v_\pm^t(q, p_0))
\end{align*}
we have
\begin{equation*}
\sup_{t \ge 0} |\tilde \varphi_\pm^t - \mathrm{Id}|_{L(\T^n \times B_{1-\delta})} = \sup_{t \le 0} | \varphi_\pm^t - \mathrm{Id}|_{L(\T^n \times D')} < C_0 \varepsilon.
\end{equation*}
We recall that $\mathrm{Leb}(B_1\backslash D) < \mu$,
\begin{eqnarray}
\label{primafurmulettamisuradelcazzo}
\mathrm{Leb}\left(\tilde \varphi^0_\pm(\T^n \times B_{1-\delta})\backslash \varphi^0_\pm(\T^n \times D')\right) &=& \mathrm{Leb}\left(\tilde \varphi^0_\pm(\T^n \times B_{1-\delta})\backslash\tilde \varphi^0_\pm(\T^n \times D')\right) \nonumber \\
&\le& C_\varepsilon \mathrm{Leb}\left(B_{1-\delta}\backslash D'\right) \nonumber\\
&=& C_\varepsilon \mathrm{Leb}\left(B_{1-\delta}\backslash \left(B_{1-\delta} \cap D\right)\right)\nonumber\\
&\le& C_\varepsilon \mathrm{Leb}\left(B_1\backslash D\right)< C_\varepsilon \mu
\end{eqnarray}
with a constant $C_\varepsilon$ converging to $1$ if $\varepsilon \to 0$. 

We observe that 
\begin{eqnarray}
\label{furmulettamisuradelcazzo}
\left(\T^n \times B_1 \right)\backslash \varphi^0_\pm(\T^n \times D') &=& \left(\T^n \times B_1 \right) \backslash \tilde \varphi^0_\pm(\T^n \times B_{1-\delta}) \nonumber \\
&+& \tilde \varphi^0_\pm(\T^n \times B_{1-\delta}) \backslash \varphi^0_\pm(\T^n \times D').
\end{eqnarray}
Thanks to Proposition \ref{propdiAbedLip} and the special form of $\tilde \varphi_\pm^0$
\begin{equation*}
\T^n \times B_{1-2\delta} \subset \tilde \varphi_\pm^0\left(\T^n \times B_{1-\delta}\right).
\end{equation*}
Then, by the latter,~\eqref{primafurmulettamisuradelcazzo} and~\eqref{furmulettamisuradelcazzo}
\begin{eqnarray}
\label{furmulettamisuradelcazzo2}
\mathrm{Leb}\left(\left(\T^n \times B_1 \right)\backslash \varphi^0_\pm(\T^n \times D')\right) &=& \mathrm{Leb}\left(\left(\T^n \times B_1\right) \backslash \tilde \varphi^0_\pm(\T^n \times B_{1-\delta}) \right) + C_\varepsilon \mu \nonumber\\
&\le& \mathrm{Leb}\left(\left(\T^n \times B_1 \right) \backslash (\T^n \times B_{1-2\delta}) \right) + C_\varepsilon \mu \nonumber\\
&\le& 2\mu
\end{eqnarray}
for $\varepsilon$ sufficiently small. Now, let us introduce the following set
\begin{equation*}
\mathcal{W} =  \varphi^0_+(\T^n \times D')\cap\varphi^0_-(\T^n \times D')
\end{equation*}
and by~\eqref{furmulettamisuradelcazzo2}
\begin{eqnarray*}
\mathrm{Leb}\left(\left(\T^n \times B_1\right)\backslash \mathcal{W}\right)&\le& \mathrm{Leb}\left(\left(\T^n \times B_1\right)\backslash  \varphi^0_+(\T^n \times D')\right)\\
&+& \mathrm{Leb}\left(\left(\T^n \times B_1\right)\backslash \varphi^0_-(\T^n \times D')\right) \le 4\mu.
\end{eqnarray*}
This concludes the proof of this theorem because, for all $(q, p_0) \in \mathcal{W} = \varphi^0_+(\T^n \times D')\cap \varphi^0_-(\T^n \times D')$, there exist $(q_+, p_{0+})$, $(q_-, p_{0-}) \in \T^n \times D'$ such that 
\begin{equation*}
\varphi_+^0(q_+, p_{0+}) = (q,p_0) = \varphi_-^0(q_-, p_{0-}).
\end{equation*}
Therefore, by Proposition \ref{proplink2BM}, there exists a biasymptotically quasiperiodic solution $g(t)$ associated to $(X_H, X_{\tilde h}, \varphi_{0, p_{0\pm}})$ such that $g(0) = (q,p_0)$.

\section{Proof of Theorem \ref{Thm4BM}}\label{ProofBM4}

The proof of this theorem is the same as Theorem \ref{Thm3BM}. However, we have some obvious differences in the estimation of the solution of the homological equation. Here, we prove the existence of a family of positive $C^\sigma$-asymptotic KAM tori $\psi^t_+$ parametrized by $p_0 \in D'$. Similarly, we have the claim concerning the existence of a family of negative $C^\sigma$-asymptotic KAM tori. In what follows, we drop the subscript $+$ to obtain a more elegant form. 

We are looking for $u$, $v : \T^n \times D' \times \R^+ \to \R^n$ such that letting
\begin{equation*}
\psi(\theta, t;p_0) = (\theta + u(\theta, t;p_0), v(\theta, t;p_0)),
\end{equation*}
for all fixed $p_0 \in D'$, $\psi$, $u$ and $v$ satisfy the following conditions
\begin{align}
\label{XHOPBM1NotInt}
&X_H(\psi(\theta, t;p_0), t;p_0) - \partial_\theta \psi(\theta, t;p_0) \omega(p_0) - \partial_t \psi(\theta, t;p_0)=0 \\
\label{XHOPBM2NotInt}
& \lim_{t \to +\infty}|u_{p_0}^t|_{C^\sigma}=0, \quad \lim_{t \to +\infty}|v_{p_0}^t|_{C^\sigma}=0, 
\end{align}
for all $(\theta, t) \in \T^n \times \R^+$.
To this end, given $\sigma \ge 1$, let us introduce the following definitions
\begin{definition}
Let $\mathcal{D}^+_{\sigma}$ be the space of functions $f$ defined on $\T^n  \times \R^+ \times D'$ such that $f \in C(\T^n  \times \R^+ \times D')$ and $f_p^t \in C^\sigma(\T^n)$ for all $(t, p) \in \R^+ \times D'$.
\end{definition}
For all $f \in \mathcal{D}^+_{\sigma}$ and $l >1$, we define the following norms
\begin{align*}
&|f|_{\sigma , l, L(D')} = \sup_{(t, p) \in \R^+ \times D'}|f^t_p|_{C^\sigma}(1 + |t|^l) + \sup_{t \in  \R^+}|f^t|_{L(D')}(1 + |t|^{l-1}),\\
&|f|_{\sigma , 0, L(D')} = \sup_{(t, p) \in \R^+ \times D'}|f^t_p|_{C^\sigma}+ \sup_{t \in  \R^+}|f^t|_{L(D')}.
\end{align*}
These norms satisfy the properties in Proposition \ref{normpropertiesBMNotInt}. As one can expect, we define the following subset of $\mathcal{D}^+_{\sigma}$

\begin{definition}
Given $\sigma \ge 1$ and an integer $k \ge 0$, we define $\mathcal{\bar D}^+_{\sigma, k}$ the space of functions $f$ such that 
\begin{equation*}
f \in \mathcal{D}^+_{\sigma +k}, \hspace{2mm} \mbox{and} \hspace{2mm} \partial^i_q f \in \mathcal{D}^+_{\sigma + k -i}
\end{equation*}
for all $0 \le i \le k$. 
\end{definition} 
Moreover, for all $f \in \mathcal{\bar D}^+_{\sigma, k}$ and $l>1$, we consider
\begin{equation*}
\left \|f \right \|^+_{\sigma, k, l, L(D')} = \max_{0 \le i \le k}|\partial_q^i f|^+_{\sigma+k-i, l, L(D')}.
\end{equation*} 

Now, let $\sigma \ge 1$ and $l >1$ be the positive parameters introduced by~\eqref{H4BM}. We define the following Banach spaces $\left(\mathcal{A}, |\cdot |\right)$, $\left(\mathcal{B}, |\cdot |\right)$, $\left(\mathcal{U}, |\cdot |\right)$, $\left(\mathcal{V}, |\cdot |\right)$, $\left(\mathcal{Z}, |\cdot |\right)$ and $\left(\mathcal{G}, |\cdot |\right)$ 

\begin{eqnarray*}
\mathcal{A} &=& \Big\{a : \T^n \times \R^+ \times D' \to \R  \hspace{1mm}| \hspace{1mm} a \in \mathcal{\bar D}^+_{\sigma, 2} \\
&& \mbox{and} \hspace{1mm} |a| = |a|^+_{\sigma+2,0, L(D')} + \left \|\partial_\theta a\right\|^+_{\sigma,1, l+2, L(D')} < \infty\Big\}\\
\mathcal{B} &=& \Big\{b : \T^n \times \R^+ \times D'  \to \R^n  \hspace{1mm}| \hspace{1mm} b \in \mathcal{\bar D}^+_{\sigma, 2,}, \hspace{1mm} \mbox{and} \hspace{1mm} |b| = \left \|b\right\|^+_{\sigma,2, l+1, L(D')} < \infty\Big\}\\
\mathcal{U} &=& \Big\{u :\T^n \times \R^+ \times D'  \to \R^n  \hspace{1mm}| \hspace{1mm} u ,  \left( \nabla_{\theta t} u\right) \Omega \in \mathcal{D}^+_{\sigma}  \\
&& \mbox{and} \hspace{1mm} |u| =\max \{|u|^+_{\sigma,l, L(D')},| \left( \nabla_{\theta t} u \right) \Omega|^+_{\sigma,l+1, L(D')} \} < \infty\Big\}\\
\mathcal{V} &=& \Big\{v : \T^n \times \R^+ \times D'  \to \R^n  \hspace{1mm}| \hspace{1mm} v,  \left( \nabla_{\theta t} v\right) \Omega \in \mathcal{D}^+_{\sigma}  \\
&& \mbox{and} \hspace{1mm} |v| =\max \{|v|^+_{\sigma,l+1, L(D')},| \left( \nabla_{\theta t} v \right) \Omega|^+_{\sigma,l+2, L(D')} \} < \infty\Big\}\\
\mathcal{Z} &=& \Big\{z : \T^n \times \R^+ \times D'  \to \R^n  \hspace{1mm}| \hspace{1mm} z \in \mathcal{D}^+_{\sigma},  \hspace{1mm} \mbox{and} \hspace{1mm} |z| =|z|^+_{\sigma,l+1, L(D')} < \infty\Big\}\\
\mathcal{G} &=& \Big\{g : \T^n \times \R^+ \times D'  \to \R^n  \hspace{1mm}| \hspace{1mm} g \in \mathcal{D}^+_{\sigma},  \hspace{1mm} \mbox{and} \hspace{1mm} |g| =|g|^+_{\sigma,l+2, L(D')} < \infty\Big\}
\end{eqnarray*} 

 Let $M_n$ be the set of the $n$-dimensional matrices and $\Upsilon \ge 1$ the positive parameter in~\eqref{H4BM}. We introduce another Banach space $\left(\mathcal{M}, |\cdot |\right)$, such that 
\begin{eqnarray*}
\mathcal{M} &=& \Big\{m : \T^n \times B_\delta \times \R^+ \times D' \to M_n  \hspace{1mm}| \hspace{1mm}  \partial^i_{\theta I}m \in C(\T^n \times B_\delta \times \R^+ \times D') \hspace{1mm}\\
&&\mbox{for all $0 \le i \le 2$,}  \hspace {2mm}m^t_{p_0} \in C^{\sigma+2}( \T^n \times B_\delta) \hspace{1mm} \mbox{for all fixed $(t, p_0) \in \R^+ \times D'$} \\
&&\mbox{and}\hspace{1mm} |m| = \hspace{1mm}\sup_{(t, p_0) \in \R^+ \times D'}|m^t_{p_0}|_{C^{\sigma + 2}} + \sup_{0 \le i \le 2}\left(\sup_{t \in \R^+} \left|\partial^i_{\theta I}m^t\right|_{L(D')}\right) \le 2\Upsilon\Big\}.
\end{eqnarray*}

Let $\mathcal{F}$ be the following functional 
\begin{equation*}
\mathcal{F} :  \mathcal{A} \times \mathcal{B} \times \mathcal{M} \times \mathcal{M} \times  \mathcal{U} \times  \mathcal{V} \longrightarrow \mathcal{Z} \times \mathcal{G}\\
\end{equation*}
\begin{equation*}
\mathcal{F}(a, b,m , \bar m,  u, v) = (F_1(b, \bar m, u,v), F_2(a,b,m,u,v))
\end{equation*}
with 
\begin{eqnarray*}
F_1(b,\bar m ,u,v) &=& b \circ \tilde u + \left(\bar m \circ \tilde \psi \right) v - \left(\nabla_{\theta t} u \right)\Omega,\\
F_2(a,b, m,u,v) &=& \partial_\theta a\circ \tilde u  + \left(\partial_\theta b \circ \tilde u \right) v  + \left(\partial_\theta m \circ \tilde \psi \right) \cdot v^2 + \left(\nabla_{\theta t} v \right)\Omega .
\end{eqnarray*}
It is obtained by~\eqref{XHOPBM1NotInt}. For fixed $m$, $\bar m \in \mathcal{M}$ and for $(a,b)$ sufficiently close to $(0,0)$, we are looking for $(u, v) \in \mathcal{U} \times \mathcal{V}$ satisfying~\eqref{XHOPBM2NotInt} such that $\mathcal{F} (a,b,m, \bar m,u,v) = 0$.
Following the lines of the proof of Theorem \ref{Thm3BM}, one can prove that $\mathcal{F}$ is well-defined, continuous, differentiable with respect to the variables $(u,v)$ with $D_{(u,v)} \mathcal{F}$ continuous. Moreover, for all fixed $m$, $\bar m \in \mathcal{M}$
\begin{equation*}
D_{(u,v)}  \mathcal{F}(0,0,m,\bar m,0,0)(\hat u, \hat v)  = (\bar m_0 \hat v - (\nabla_{\theta t} \hat u) \Omega, (\nabla_{\theta t} \hat v) \Omega).
\end{equation*}
is invertible. The proof relies on the solution of the following homological equation.
Given $\sigma \ge 1$ and $l>1$, we consider the following equation in the unknown $\varkappa : \T^n \times \R^+ \times D' \to \R$
\begin{equation}
\label{HEBMNotInt}
\begin{cases}
 \omega(p_0) \cdot \partial_q \varkappa(\theta, t ; p_0) + \partial_t \varkappa(\theta, t ; p_0) = g(\theta, t ; p_0),\\
g \in \mathcal{D}^+_{\sigma}, \quad |g|^+_{\sigma, l + 1, L(D')} < \infty,\\
\omega : D' \longrightarrow \R^n,\\
\omega \in C(D'), \quad |\omega|_{L(D')} < \infty.
\end{cases}
\tag{$HE_B$}
\end{equation}

\begin{lemma}[\textbf{Homological Equation}]
\label{homoeqlemmaBMNotInt} 
There exists a unique solution $\varkappa \in \mathcal{D}^+_{\sigma}$ of~\eqref{HEBMNotInt} such that, for all fixed $p_0 \in D'$,
\begin{equation}
\label{asymptoticBMNotInt}
\lim_{t \to +\infty} |\varkappa^t_{p_0}|_{C^0} = 0.
\end{equation}
Moreover, 
\begin{equation*}
|\varkappa|^+_{\sigma, l, L(D')} \le C(l, |\omega|_{L(D')}) |g|^+_{\sigma, l+1, L(D')}
\end{equation*}
for a suitable constant $C(l, |\omega|_{L(D')})$ depending on $n$, $l$ and $|\omega|_{L(D')}$. 
\end{lemma}
\begin{proof}
By Lemma \ref{homoeqlemmaBM}, we know that
\begin{equation*}
\label{chiNotInt}
\varkappa(\theta, t ; p_0) =-\int_t^{+\infty} g(\theta+ \omega(p_0)(\tau - t),\tau ;p_0)d\tau
\end{equation*}
is the unique solution of (\ref{HEBMNotInt}) satisfying~\eqref{asymptoticBMNotInt}. Concerning the estimates, similarly to the proof of Lemma \ref{homoeqlemmaBM}, we have
\begin{equation*}
\sup_{(t, p_0) \in  \R^+ \times D'}|\varkappa_{p_0}^t|_{C^\sigma}(1 + t^l) \le C(l)|g|^+_{\sigma,l+1, L(D')}.
\end{equation*}
It remains to estimate the second member of the norm.  By Lemma \ref{lemmastimeBMjo}, for all $(\theta, t,p_{0_1})$, $(\theta, t,p_{0_2}) \in \T^n \times \R^+ \times D'$ with $p_{0_1} \ne p_{0_2}$, 
\begin{align*}
&{|\varkappa(\theta, t; p_{0_1}) -\varkappa(\theta, t; p_{0_2})| \over |p_{0_1} - p_{0_2}|} \\
&\le \int_t^{+\infty}{|g(\theta + \omega(p_{0_2})(\tau - t), t, p_{0_2}) -g(\theta + \omega(p_{0_2})(\tau - t), t, p_{0_1})| \over |p_{0_1} - p_{0_2}|}d\tau \\
&\le \int_t^{+\infty}{|g(\theta + \omega(p_{0_2})(\tau - t), t, p_{0_1}) -g(\theta + \omega(p_{0_1})(\tau - t), t, p_{0_1})| \over |p_{0_1} - p_{0_2}|}d\tau \\
&\le C\left(\sup_{t \in \R^+}|g^t|_{L(D')}\left(1 + t^l\right)\right)\int_t^{+\infty}{1 \over 1 + \tau^l}d\tau\\
&+C\left(\sup_{(t,p_0) \in \R^+ \times D'}|g^t_p|_{C^1}\left(1 + t^{l+1}\right)\right)|\omega|_{C^1}\int_t^{+\infty}{\tau-t \over 1 + \tau^{l+1}}d\tau,\\
&\le C(l, |\omega|_{C^1}) {|g|^+_{1, l+1, L(D)} \over 1 + t^{l-1}}.
\end{align*}
Taking the sup for all $\theta \in \T^n$, $p_{0_1}$, $p_{0_2} \in D'$ with $p_{0_1} \ne p_{0_2}$, and then for all $t \in \R^+$ on the left-hand side of the latter, we conclude the proof of this lemma.
\end{proof}

We proved that the functional $\mathcal{F}$ satisfies the hypotheses of the implicit function theorem. Therefore, following the lines of the proof of Theorem \ref{Thm3BM}, one can conclude the proof of Theorem \ref{Thm4BM}.

\section*{Acknowledgement}

\textit{These results are part of my PhD thesis, which I prepared at Université Paris-Dauphine. I want to thank my thesis advisors, Abed Bounemoura and Jacques Féjoz. Without their advice and support, this work would not exist.}

\textit{This project has received funding from the European Union’s Horizon 2020 research and innovation programme under the Marie Skłodowska-Curie grant agreement No 754362}. \includegraphics[scale=0.01]{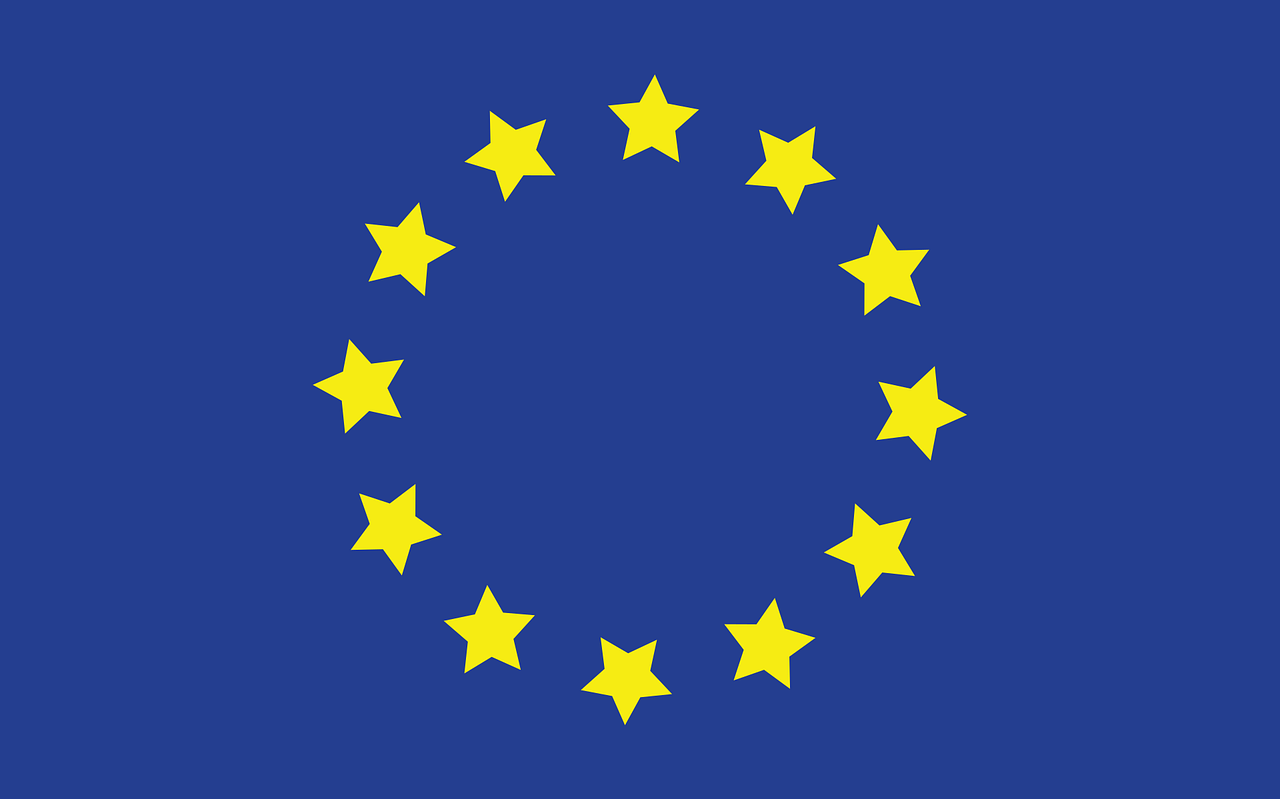}

\bibliographystyle{amsalpha}
\bibliography{Biasymptoticfinal}

\end{document}